\documentclass[10pt]{article}
\font\smallit=cmti10

\usepackage{lmodern}
\usepackage[T1]{fontenc}
\usepackage{xcolor}
\usepackage[export]{adjustbox}
\usepackage{forest}
\usepackage{amssymb,amsmath,amsthm,enumitem} 
\usepackage{subcaption,booktabs}
\usepackage{amstext}
\usepackage{amsmath}
\usepackage{amsfonts}
\usepackage{array}
\usepackage{url}
\usepackage{indentfirst}
\usepackage{tikz}
\usepackage{filecontents}
\usepackage[margin=1.5in]{geometry}    
\usepackage{graphicx}
\usepackage[utf8]{inputenc}
\usepackage{hyperref}
\usetikzlibrary{shapes.multipart}
\tikzset{block/.style={
        font=\sffamily,
        draw=black,
        thin,
        fill=white!50,
        rectangle split,
        rectangle split horizontal,
        rectangle split parts=#1,
        outer sep=0pt},
        }

\usepackage{enumerate} 
\usepackage[justification=centering]{caption} 
\usepackage{caption}
\usetikzlibrary{arrows.meta}

\makeatletter

\renewcommand\section{\@startsection {section}{1}{\z@}
{-30pt \@plus -1ex \@minus -.2ex}
{2.3ex \@plus.2ex}
{\normalfont\normalsize\bfseries\boldmath}}

\renewcommand\subsection{\@startsection{subsection}{2}{\z@}
{-3.25ex\@plus -1ex \@minus -.2ex}
{1.5ex \@plus .2ex}
{\normalfont\normalsize\bfseries\boldmath}}

\renewcommand{\@seccntformat}[1]{\csname the#1\endcsname. }

\makeatother

\newtheorem{theorem}{Theorem}

\newtheorem{lemma}{Lemma}
\newtheorem{observation}{Observation}   
\newtheorem{corollary}{Corollary} 

\newtheorem{proposition}{Proposition}

\theoremstyle{definition}

\newtheorem{remark}{Remark}

\def\ZZ{\mathbb{Z}}

\title{On pairs of triangular numbers 
whose product is a perfect square 
and pairs of intervals of successive integers with equal sums of squares}
\begin{document}
\maketitle

\vskip 15pt
\begin{center}
{\bf Vladimir Gurvich}\\
{\smallit National Research University Higher School of Economics 
(HSE), Moscow, Russia, 
and RUTCOR, Rutgers,the State University of NJ, USA}\\
{\tt vgurvich@hse.ru}, {\tt vladimir.gurvich@gmail.com}\\
\vskip 10pt

{\bf Mariya Naumova}\\
{\smallit Rutgers Business School, Rutgers University, Piscataway, NJ, United States}\\
{\tt mnaumova@business.rutgers.edu}\\
\end{center}

\begin{abstract} 
In 1778 Leonhard Euler characterized triangular numbers that are perfect squares. 
Obviously, the product of any two such numbers is a perfect square too. 
Yet, there are many other solutions, that is, pairs $(k,k')$ 
such that  $k(k+1)k'(k'+1)$ is a perfect square. 
We give explicit formulas characterizing all these square triangular pairs 
by means of some integer positive polynomials, which is the primary novelty of our work. 

This result allows us to find all pairs of intervals of successive integers 
with equal sums of squares in case when the lengths of two intervals in a pair differ by 1. 
It is known that there is a one-to-one correspondence between 
the square triangular numbers and nearly isosceles Pythagorean triples:  
$n^2 + (n+1)^2 = N^2$. 
Both are generated by the same Fermat-Pell recursion. 
This is a special case of our result, when the lengths of the two intervals are 2 and 1. 
\newline 
SUBJECT CLASSIFICATION: number theory, 
triangular square numbers,  
(near-isosceles) Pythagorean triples, Pell's equation.
\end{abstract} 

\section{Introduction}
Denote by $\ZZ_+$ the set of positive integer numbers; 
a number $n \in \ZZ_+$  is called:  
\begin{itemize}
\item [] {\em triangular}  if  $n = k(k + 1)/2$ for some $k \in \ZZ_+$;
\item [] {\em square} or {\em a perfect square} if $n = k^2$ for some $k \in \ZZ_+$;   
\item [] {\em square-free} if $k^2$ is a divisor of $n$ only for $k=1$; 
in other words, if $n$ is the product of pairwise distinct primes.
\end{itemize}

By convention, every prime number is square-free and 1 is square-free too, 
it corresponds to the empty set of primes. 
Thus, 1 is the only number that is square and square-free simultaneously.
Obviously, any $n \in \ZZ_+$  can be uniquely represented 
as a product $n = n' \times n''$, 
where  $n'$ is square and $n''$  is square-free. 
    
The square triangular numbers were characterized by Leonhard Euler~ 
\cite[Exemplum~I]{Eul1778}.  
It is also known that 
these numbers are in one-to-one correspondence with the so-called 
{\em near-isosceles Pythagorean triples} $(k,k+1,n)$  such that  $k^2 + (k+1)^2 = n^2$. 

A {\em quadratic number} is 
the product  $\Pi(k,k') = k(k+1)k'(k'+1)$  
for some  $k, k' \in \ZZ_+$, which may coincide. 
In other words, it is the product of two triangular numbers and 4.  
A pair $(k,k')$ is called {\em square} if $\Pi(k,k')$  is square. 
Clearly, $(k,k')$ is square if $k = k'$  or when both triangular numbers 
$k(k+1)/2$ and $k'(k'+1)/2$ are square. 
However, there are many other solutions. 

In this article we characterize all square pairs. 
We will partition $\ZZ_+$ into an infinite family 
of infinite subsets such that 
$\Pi(k,k')$ is square if and only if 
$k$ and $k'$ belong to the same subset.
For example, all triangular square  numbers form one such subset. 

\medskip 

Every $n  \in \ZZ_+$ admits a unique prime factorization. 
The product of all primes that appear 
in this factorization an odd number of times 
(in odd powers) will be called the {\em odd radical} of $n$ and denote $oddrad(n)$.  
By convention, $oddrad(n) = 1$  if  $n$ is a perfect square, 
otherwise  $oddrad(n) > 1$. 
By definition, every odd radical is square-free. 

\begin{observation}
\label{o-0}    
The product $n = n' n''$  is square if and only if  $oddrad(n') = oddrad(n'')$. 
\end{observation}

\proof 
Indeed, a number $n$ is square if and only if 
in its factorization 
every prime appears in an even power; 
in other words, if $oddrad(n) = 1$.
\qed 

\begin{corollary}
If $n n'$ and $n n''$ are both  square   
then $n' n''$  is square, too.
\qed 
\end{corollary}

In other words, the relation is transitive. 
It is also idempotent, since  
$n^2$ is square for every  $n$. This is a tautology. 
Thus, $\ZZ_+$ is partitioned into equivalence classes such that 
all numbers of a class have the same odd radical.  
Then, the product of two numbers is square if and only if 
they belong to one class. 

\begin{observation}
\label{o-1}    
Each class is infinite and there are infinitely many classes. 
\end{observation}

\proof 
There are infinitely many odd radicals 
and each one can be multiplied 
by infinitely many perfect squares.
\qed 

\medskip 

We will restrict the above simple theory 
from $\ZZ_+$  to the set $T$ of triangular numbers and  
show that Observations \ref{o-0} and \ref{o-1} still hold.  
More precisely, we partition $T$  into subsets such that 
all numbers from one  subset have the same odd radical. 
Then, $\Pi(k,k') = k(k+1)k'(k'+1)$ is square if and only if 
the corresponding two triangular numbers $k(k+1)/2$ and $k'(k'+1)/2$  are from the same subset. 
We will show that each subset is infinite and there are infinitely many of them. 
For example, all square triangular numbers form one subset with $oddrad(n) = 1$, 
since each prime appears in an even power in the factorization of any number of this subset. 

Any triangular number can be uniquely represented as 

\begin{equation}
\label{eq-a}
k(k+1)/2 = d t^2, \text{ where } d = oddrad(k(k+1)/2) \;\; \text{ is square-free }.
\end{equation}

Multiplying by 8 we transform this equation to  

\begin{equation}
\label{eq-b}
(2k + 1)^2 - 2d (2t)^2 = 1.
\end{equation}

Thus, we obtain a special case of the well-known Fermat-Pell equation

\begin{equation}
\label{eq-Pell}
x^2 - D y^2 = 1,
\end{equation}

\noindent 
with  $x = 2k+1, y=2t$, and $D = 2d$.   
Fermat-Pell equation has a very long history and is very well studied.
Obviously, it has no solution if $D$ is square. 
Otherwise, it has infinitely many solutions  $x_i, y_i; i = 0, 1, 2, \dots$. 
Numbers $x_i$ and $y_i$ increase monotonically (and very rapidly) as $i$ increases. 
They are efficiently characterized in several ways: 
recursively, by continued fractions, and  by explicit formulas. 
We will not repeat these classical results here and refer the reader to 
\cite{Dic05,HM85,Mah84,Mar50,McD96,Wei06,Wiki-Pell}.
Equation (\ref{eq-a}) was solved by Euler; see Section \ref{s-2}. 

Obviously, $D = 2d$  is square if and only if $d=2$,  
since $d$ is square-free. 
Thus, we obtain

\begin{proposition}
Equation (\ref{eq-b}) has no solutions if and only if $d=2$.  
In other words, each square-free number, except 2, is 
an odd radical of a triangular number. 
If $d \neq 2$ then (\ref{eq-b}) has infinitely many solutions. 
\qed 
\end{proposition}


Given a square-free  $d$  distinct from 2, 
we can find the set of all triangular numbers with the odd radical $d$. 
To do so, we solve (\ref{eq-a})  using the Fermat-Pell theory. 
Any two solutions (not necessarily distinct) provide a quadratic square number. 

In this article we apply a different approach. 
Given an arbitrary  $k = k_0 \in \ZZ_+$, 
we provide explicit formulas for all $k_i$  such that 
$oddrad(k_i(k_i+1)/2) = oddrad(k_0(k_0+1)/2)$. 
Then, quadratic number $\Pi(k_i,k_j) = k_i(k_i+1)k_j(k_j+1)$ is a perfect square  
for all $i,j = 0,1,2,...$  
Here $k_i = P_i(k)$  given by a polynomial of degree $i+1$. 
For $i = 0,1,\dots,6$ we obtain:   
$$k_0 = k, \; k_1 = (2k+1)^2-1, \; k_2 = (k + 1)(4k + 1)^2 - 1 , \;k_3 = (8k^2+8k+1)^2 - 1,$$    
$$k_4 = (k + 1)((4k)^2 + 3(4k) + 1)^2 - 1, \;   
k_5 = (2k + 1)^2(16k^2 + 16k + 1)^2 - 1, \dots$$
$$k_6 = (k + 1)((4k)^3 + 5(4k)^2 + 6(4k) + 1)^2 - 1,\dots$$ 

We derive a general formula working for all  $i$ as follows. 
Given $k$ and  $i$, define  $k_i = P_i(k)$  by the next equations 
(distinct for odd and even $i$) 
of degree $i = 2\ell$ or $2\ell + 1$ 
with  coefficients 
$(a_t \mid t = 0, \ldots, \ell; \;  b_u \mid  u = 0, \ldots, \ell+1)$: 

\medskip 

$k_{2 \ell} + 1 = 
k [a_\ell k^\ell + \dots + a_1 k + a_0]^2 +1 = 
(k+1) [b_\ell k^\ell + \dots + b_1 k + b_0]^2$  and

\medskip 

$k_{2\ell+1} + 1 = 
 k(k+1) [a_\ell k^\ell + \dots + a_1 k + a_0]^2 + 1 = 
 [b_{\ell+1} k^{\ell+1}+b_\ell k^\ell + \dots + b_1 k + b_0]^2$.  

\medskip 

We prove that for every $i$ and $k$ there exist unique integer positive 
$a,b$, and  $k_i$  such that $(k, k_i)$ is a square pair, 
that is,  $\Pi(k, k_i) = k (k+1) k_i(k_i+1)$  is a perfect square. 
Numbers $k_i = k_i(k)$ are given in the Table \ref{t2} for $i \leq 8$ and $k \leq 8$. 

\begin{table}[htbp]
\footnotesize
\setlength{\tabcolsep}{1.2pt}
\begin{center}
\begin{tabular}{p{0.1cm} c c c c c c c c c}
$k$ & $k_1$ & $k_2$ & $k_3$ & $k_4$ & $k_5$ & $k_6$ &  $k_7$ & $k_8$ \\
\hline
$1$ & $8$ & $49$ & $288$ & $1,681$ & $9,800$ & $57,121$ & $332,928$ & $1,940,449$\\ 
$2$ & $24$ & $242$ & $2,400$ & $23,762$ & $235,224$ & $2,328,482$ & $23,049,600$ & $228,167,522$\\ 
$3$ & $48$ & $675$ & $9,408$ & $131,043$ & $1,825,200$ & $25,421,763$ & $354,079,488$ & $4,931,691,075$\\ 
$4$ & $80$ & $1,444$ & $25,920$ & $465,124$ & $8,346,320$ & $149,768,644$ & $2,687,489,280$ & $48,225,038,404$\\
$5$ & $120$ & $2,645$ & $58,080$ & $1,275,125$ & $27,994,680$ & $614,607,845$ & $13,493,377,920$ & $296,239,706,405$\\
$6$ & $168$ & $4,374$ & $113,568$ & $2,948,406$ & $76,545,000$ & $1,987,221,606$ & $51,591,216,768$ & $1,339,384,414,374$\\
$7$ & $224$ & $6,727$ & $201,600$ & $6,041,287$ & $181,037,024$ & $5,425,069,447$ & $162,571,046,400$ & $4,871,706,322,567$\\
$8$ & $288$ & $9,800$ & $332,928$ & $11,309,768$ & $384,199,200$ & $13,051,463,048$ & $443,365,544,448$ & $15,061,377,048,200$\\
\hline
\end{tabular}
\end{center}
\caption{All perfect square quadratic pairs $(k, k_i)$ 
for $k \leq 8$ and $i \leq 8$. 
\newline
Each quadratic number $\Pi(k_i,k_j) = k_i(k_i+1)k_j(k_j+1); i,j \; \geq 0$, is a perfect square.} 
\label{t2}
\end{table}

The first row $k=1$  contains all numbers $k_i$ 
such that $k_i(k_i+1)/2$  is a perfect square: 
$1, 8, 49, 288, \dots$  for  $i = 0,1,2,3,\dots$
Note that the set of numbers in the row $k=8$ is a proper subset of this set.
Hence, row $k=8$ is not needed for constructing equivalence classes;  
we can skip it using the row $k=1$ instead.

In general, if a number $k_i$ appears in row $k$ 
then the set of numbers of row $k_i$ is 
a subset of the set of numbers of row  $k$.  
Our polynomial equations and formulas hold for such rows as well.  
Yet, row $k$ is essential only if $k$ did not appear 
in a previous row  $k' < k$. 
We will call such a row and number $k$ {\em basic}. 
The non-basic rows are irrelevant.  

By construction, all triangular numbers in a row have the same odd radical. 
Conversely, each square-free number, except 2, appears as 
$oddrad(k(k+1)/2$  for some basic $k$.
Thus, there is a one-to-one correspondence between 
the basic rows and square-free numbers, except  2. 
In other words, our construction gives 
all pairs of triangular numbers whose product is a perfect square. 


\medskip 

In the last section we consider an important application of this result 
and find all pairs of intervals of successive integers 
with equal sums of squares provided 
the lengths of two intervals in a pair differ by 1. 

\section{Square triangular numbers} 
\label{s-2}
Here we characterize square triangular numbers. 
In Section \ref{near-iso} we show that 
there is a one-to-one correspondence between them 
and the near-isosceles Pythagorean triples, $n^2 + (n+1)^2 = N^2$. 
In both cases the characterization is based on 
the same Fermat-Pell equation; 
see \cite{Bar03,BT01,Fin87,FL68,Hin83,PSJW62} for more details. 

As we mentioned, the square triangular numbers 
were characterized by Euler in 1778 \cite{Eul1778}, 
while the near-isosceles Pythagorean triples appear  
in the ancient Babylonian clay tablet known as Plimpton 322, dated to around 1800 BC, 
that is, long before Pythagoras \cite{MW17}.  
Yet, they were characterized much later,  
as far as we know, in Beiler's book \cite{Bei64}; see also \cite{Her83}.

\subsection{Recursive and explicit characterizations} 
Consider the following ``Fibonacci-like'' sequence 
$$(x_0,x_1, x_2, \dots) = 
(0, 1, 2, 5, 12, 29, 70, 169, 408, 985, \; 2,378, \; 5,741, \dots)$$   
initialized by 
$x_0 =  0, x_1 = 1$  and defined by recursion 
$x_{i+2} = 2 x_{i+1} + x_i$  for $i = 1, 2, \ldots$. Its characteristic equation 
$x^2 - 2x - 1 = 0$
with discriminant $D=8$ and two roots $x = 1 \pm \sqrt{2}$    
provides the following explicit formula: 
\begin{equation}
\label{x_i}
x_i = \frac{\sqrt{2}}{4}\left((1+\sqrt{2})^i - (1-\sqrt{2})^i\right), i = 0, 1, \ldots
\end{equation}

It is obvious and well known \cite{Mar50} that 
$x_{i+1}/x_i \rightarrow 1 + \sqrt{2}$ is as $i \rightarrow \infty$.

The above recursion is a special case of the Fermat-Pell recurrence \cite{HM85, Mah84, Wei06}.

\begin{lemma} 
\label{t-tr-sq}
For $i=0, 1, \ldots$, if $t_i = 2x_i^2$ when $i$ is even and $t_i = (x_{i-1}+x_{i})^2$ 
when $i$ is odd, then triangular number $t_i(t_i+1)/2$ is a perfect square.
\end{lemma} 

\proof 
Assume $i$ is even. Then
$t_i(t_i+1)/2 = x_i^2(2x_i^2+1)$, and by (\ref{x_i}),
\begin{equation}
\begin{split}
t_i(t_i+1)/2 &= \frac{1}{8}\left((1+\sqrt{2})^i - (1-\sqrt{2})^i\right)^2\left(2\cdot \frac{1}{8}\left((1+\sqrt{2})^i-(1-\sqrt{2})^i\right)^2+1\right)\\
&= \frac{1}{8}\left((1+\sqrt{2})^i - (1-\sqrt{2})^i\right)^2\cdot\frac{1}{4}\left((1+\sqrt{2})^i+(1-\sqrt{2})^i\right)^2\\
&= \frac{1}{32}\left((1+\sqrt{2})^{2i} - (1-\sqrt{2})^{2i}\right)^2.
\end{split}\nonumber 
\end{equation}

Assume $i$ is odd. Then
$t_i(t_i+1)/2 = (x_{i-1}+x_i)^2\left((x_{i-1}+x_i)^2+1\right)/2$. By (\ref{x_i}),

\begin{equation}
\begin{split}
x_{i-1}+x_i & = \frac{\sqrt{2}}{4}\left(\frac{(1+\sqrt{2})^{i}}{1+\sqrt{2}} - \frac{(1-\sqrt{2})^i}{1-\sqrt{2}} + (1+\sqrt{2})^i - (1-\sqrt{2})^i\right)\\
& = \frac{\sqrt{2}}{4}\left(\frac{2+\sqrt{2}}{1+\sqrt{2}}(1+\sqrt{2})^{i} - \frac{2-\sqrt{2}}{1-\sqrt{2}}(1-\sqrt{2})^i\right)\\
& = \frac{\sqrt{2}}{4}\left((1+\sqrt{2})^{i} \cdot \sqrt{2}+ (1-\sqrt{2})^i\cdot \sqrt{2}\right)\\
& = \frac{1}{2}\left((1+\sqrt{2})^{i} + (1-\sqrt{2})^i\right),
\end{split}\nonumber 
\end{equation}

\begin{equation}
(x_{i-1}+x_i)^2 = \frac{1}{4}\left((1+\sqrt{2})^{2i} + (1-\sqrt{2})^{2i}-2\right),
\nonumber 
\end{equation}
\noindent
so that

\begin{equation}
t_i(t_i+1)/2 = \frac{1}{2}(x_{i-1}+x_i)^2((x_{i-1}+x_i)^2+1)
= \frac{1}{32}\left((1+\sqrt{2})^{2i} - (1-\sqrt{2})^{2i}\right)^2.
\nonumber 
\end{equation}

Thus, in both cases,

\begin{equation}
\label{sq_r}
\sqrt{t_i(t_i+1)/2} = \frac{(1+\sqrt{2})^{2i} - (1-\sqrt{2})^{2i}}{4\sqrt{2}}.   
\end{equation}
Clearly, (\ref{sq_r}) is an integer for $i=0$. 
Assuming (\ref{sq_r}) is an integer for $i=k$, 
\begin{equation}
\begin{split}
\sqrt{t_{k+1}(t_{k+1}+1)/2} & = \frac{(1+\sqrt{2})^{2k+2} - (1-\sqrt{2})^{2k+2}}{4\sqrt{2}} \\
& = \frac{(1+\sqrt{2})^{2k}(3+2\sqrt{2}) - (1-\sqrt{2})^{2k}(3-2\sqrt{2})}{4\sqrt{2}}\\
& =  \sqrt{t_k(t_k+1)/2} + \frac{1}{2\sqrt{2}}\left((1+\sqrt{2})^{2k+1} - (1-\sqrt{2})^{2k+1}\right)\\
& =  \sqrt{t_k(t_k+1)/2} + x_{2k+1},
\end{split}\nonumber 
\end{equation}
which is an integer.

\qed 

The first several examples are given in Table \ref{t1}: 

\begin{table}[h] 
\small
\begin{center}
\begin{tabular}{c c c c}
$i$ & $t_i$ & $t_i (t_i+1)/2$ & $\sqrt{t_i (t_i+1)/2}$\\ 
\hline
$0$ & $2 \cdot x_0^2 = 0$ &  $0 \cdot 1/2 = 0$& $0$\\
$1$ & $(x_0+x_1)^2 = (0+1)^2 = 1$ & $1 \cdot 2/2 = 1$ & $1$\\
$2$ & $2 \cdot x_2^2 = 2 \cdot 2^2 = 8$ &  $8 \cdot 9/2 = 36$ & $6$\\
$3$ & $(x_2+x_3)^2 = (2+5)^2 = 7^2 = 49$ & $49 \cdot 50/2 = 1,225$ & $35$\\
$4$ & $2 \cdot x_4^2 = 2 \cdot 12^2 = 288$ & $288 \cdot 289/2 = 41,616$ &  $204$\\
$5$ & $(x_4+x_5)^2 = (12+29)^2 = 41^2 = 1,681$ & $1,681 \cdot 1,682/2 = 1,413,721$ & $1,189$\\
\hline
\end{tabular}
\end{center}
\caption{First triangular numbers that are perfect squares} \label{t1}
\end{table}

\subsection{Euler's explicit solution}
Euler's solution of the equation 

\begin{equation}
t(t+1) = 2d^2.   
\label{Pell}
\end{equation}

\noindent 
is explicit and given in different form: 
\begin{equation}
    t = \frac{(3+2\sqrt{2})^n}{4}+\frac{(3-2\sqrt{2})^n}{4}-\frac{1}{2}, \hspace{5mm} d = \frac{(3+2\sqrt{2})^n}{4\sqrt{2}}-\frac{(3-2\sqrt{2})^n}{4\sqrt{2}}.\label{Eul}
\end{equation}

Euler did not prove that all solutions of \eqref{Pell} are given by this representation. 
Yet, it is not difficult to show \cite{Dic05}. It is easily seen that \eqref{Pell} is equivalent to 

\begin{equation}
(2t + 1)^2 - 2 \cdot (2d)^2 = 1,
\label{classic_pell}
\end{equation}
which is a special case of the famous Pell (Fermat-Pell) equation \cite{Wei06}. 
The solutions of the Pell's equation \eqref{classic_pell} have not only closed-form representations, 
but also recursive ones. For example, the solutions of the given equation can be written as:
\[
t_0 = 0, \quad d_0 = 0,
\]
\[
t_{k+1} = 3t_k + 4d_k + 1,
\]
\[
d_{k+1} = 2t_k + 3d_k + 1.
\]
Thus, we obtain the sequence \((t, d) = (0,0), (1,1), (8,6), (49,35), \ldots\), 
which coincides with the analytic solution.

\medskip 

It is easily proven by contradiction that 
the above construction provides all solutions.  
Given the sequences \(d_k\) and \(t_k\), suppose that there exists a solution with \(d_k < d < d_{k+1}\). 
Then, \(t_k < t < t_{k+1}\). Using the recurrence formulas in reverse, 
it is not hard to see that there must exist a solution \(\overline{d}\) 
such that \(d_{k-1} < \overline{d} < d_k\). 
Continuing with this argument, we show that there is a solution between \(d_0\) and \(d_1\), which is a contradiction.

\medskip 

The reader can find the detailed history 
of the Fermat-Pell equations in Chapter VI of Dickson's monograph \cite{Dic05}. 
Later, various explicit and recursive solutions were obtained in 
\cite{Bur10,CP95,Choi16,DN98,Hat95,Hin83,Luc61,McD96,Nyb98,PSJW62}.

\section{Equations realizing square quadratic numbers}
Obviously, the product of any two triangular square numbers is a perfect square.  
Yet, there are many other solutions. 
We will search them in the form  \(\Pi(k, k_i) = k(k+1)k_i(k_i + 1)\), 
where $k_i = P_i(k)$  are polynomials defined for even and odd $i$  by the following two identities:
\begin{equation}
\label{eq-even-new}
k_{2\ell} + 1 = 
k \left(a_\ell k^\ell + \dots + a_1 k + a_0\right)^2 + 1 = 
(k+1) \left(b_\ell k^\ell + \dots + b_1 k + b_0\right)^2,
\end{equation}
\begin{equation}
\label{eq-odd-new}
k_{2\ell+1} + 1 = 
k(k+1) \left(a_\ell k^\ell + \dots + a_1 k + a_0\right)^2 + 1 = 
\left(b_{\ell+1} k^{\ell+1} + b_\ell k^\ell + \dots + b_1 k + b_0\right)^2,
\end{equation}

\noindent
where \(\ell \geq 0\) is an integer, and \(i = 2\ell\) or \(i = 2\ell + 1\), respectively.

\medskip

We rewrite~\eqref{eq-even-new} and~\eqref{eq-odd-new} in the form:
\begin{equation}
\label{eq-ABC-even}   
k_{2\ell} + 1 = 
k \left( \sum_{t=0}^\ell a_t k^t \right)^2 + 1 = 
(k+1) \left( \sum_{t=0}^\ell b_t k^t \right)^2 = 
\sum_{t=0}^{2\ell+1} c_t k^t,
\end{equation}
\begin{equation}
\label{eq-ABC-odd}   
k_{2\ell+1} + 1 = 
k(k+1) \left( \sum_{t=0}^\ell a_t k^t \right)^2 + 1 = 
\left( \sum_{t=0}^{\ell+1} b_t k^t \right)^2 = 
\sum_{t=0}^{2(\ell+1)} c_t k^t.
\end{equation}

\begin{remark}
\label{notation}
The reader should be careful with notation. 
We use the same summation index $t$ 
for the coefficients of 
polynomials $A,B,$ and $C$ in both formulas. 
The summation is always from 0 
but the upper bound take values 
$\ell, \ell+1, 2\ell, 2\ell+1$, and $2\ell+2$ 
in (\ref{eq-ABC-even}) and (\ref{eq-ABC-odd}). 
\end{remark}

We will show that for each  $i$, even or odd, 
there exists a unique triplet of polynomials $A_i,B_i,C_i$ 
with integer positive  coefficients satisfying  
(\ref{eq-ABC-even}) or (\ref{eq-ABC-odd}), respectively, 
so that these equations can be written as:

\medskip 

$k_{2\ell} + 1  = 
k A_{2\ell}^2 +1 = (k+1) B_{2\ell}^2$ and 
$k_{2\ell + 1} + 1  = 
k(k+1) A_{2\ell+1}^2 +1 = B_{2\ell+1}^2$. 

\medskip 
By construction, the corresponding quadratic number 
is a perfect square for all $k$ and $i$: 
$$\Pi(k,k_i) = [k(k+1) k_i(k_i+1)] 
= [k(k+1)A_i B_i]^2.$$  

\begin{remark}
Note that (\ref{eq-even-new}) and (\ref{eq-odd-new})  
generate only some special perfect square quadruples. 
For example, we use identities of the form 
$$kA_i^2+1 = (k+1)B_i^2 
\;\;\; \text{or} \;\;\;
k(k+1)A_i^2+1 = B_i^2$$  
for even and odd $i$, respectively. 
One could swap them or try a more general form:  
$$kA_i^2E_i+1 = (k+1)B_i^2 E_i 
\;\; \text{or} \;\; 
k(k+1)A_i^2E_i+1 = B_i^2E_i.$$  
for some polynomial $A_i,B_i$ and $E_i$,  
which would give us perfect square quadratic numbers 
$\Pi = [k(k+1)A_iB_iE_i]^2$, too. 
However, 
such modifications do not provide new square quadratic numbers.   
Moreover, all  such numbers  
are realized by equations~\eqref{eq-even-new} and~\eqref{eq-odd-new}. 
\end{remark}

\begin{remark}
Let us note that 
\newline
$k(k+1)k(k+1) = [k(k+1)]^2$,\; 
\newline 
$k(k+1)(k+1)(k+2) = [(k+1)^2 - 1](k+1)^2$, \; 
\newline
$k(k+1)(k+2)(k+3) = (k^2 + 3k + 1)^2 - 1$. 
\newline
Thus, the first product is always perfect square, 
while the last two - never. 
\end{remark}

\section{Case of even $i$}
Given $\ell > 0, i= 2 \ell$, and $0 \leq t \leq \ell$, 
define $a_t$ and $b_t$ by the chains of identical expressions 
\begin{equation}
\label{eq-ab-even} 
a_t = 2^{2t}\frac{(\ell+t)!}{(2t+1)! (\ell-t)!} (2\ell+1),  
\;\;\; b_t = 2^{2t} \frac{(\ell+t)!}{(2t)! (\ell-t)!} = 
2^{2t} \binom{\ell+t}{2t} = 2^{2t} \binom{\ell+t}{\ell-t}.
\end{equation}

Note that $a_\ell=b_\ell=2^{2\ell}, \; 
a_0 = 2\ell+1, \; b_0 = 1$. 
Furthermore, $b_t$ is integer, by definition. 

\begin{observation}
\label{obs-integer-even}
Coefficients  $a_t$  are integer too.
\end{observation}

\proof 
We show more, namely, that $a_t /2^{2t}$  is integer. 
To do so, set 
$$d = \ell - t 
\;\; \text{and} \;\; 
M = \frac{(2t+d)!}{(2t+1)! d!} =   
\frac{(2t+2) \times \dots \times (2t+d)}
{2 \times \dots \times d}.$$
Furthermore, 
$\frac{a_t}{2^{2t}} = \frac{(2t+d)!}{(2t)! d!} \; \frac{2(t+d)+1}{2t+1} = 
M (2(t+d)+1)$.
Obviously, $a_t = 2^{2t}$ if $d=0$ and 
$a_t = 2^{2t}(2t+3)$ if $d=1$. 
If $d \geq 2$ then $M$ is a multiple of $d!$ 
unless $2t+1$  is a multiple of $d$, 
but in this case factor $2\ell+1 = 2(t+d)+1$ is a multiple of $d$ too. 
\qed    

\medskip

Introduce polynomial $C_i = \sum_{t=0}^{2\ell+1} c_t k^t$, 
with $i = 2\ell$, by the following formula: 

\begin{equation}
  c_t = 
  \begin{cases}
    \frac{2^{3t-1}} {(2t)!} 
\frac{(\ell+\frac{t}{2})!}{(\ell-\frac{t}{2})!} 
\frac{(2\ell+(t-1))!!}{(2\ell-(t-1))!!}(2\ell+1)  
& \text{for even} \;\; t \in [2;\ell],  
\;\; \text{while} \;\; c_0 = 1;\\
    \frac{2^{3t-2}} {(2t)!} 
\frac{(\ell+\frac{t-1}{2})!}{(\ell-\frac{t-1}{2})!} 
\frac{(2\ell+t)!!}{(2\ell-t)!!}(2\ell+1) 
& \text{for odd} \;\; t \in [1;\ell];\\
 2^{2t-1}\frac{(2\ell+t)!}{(2t)!(2\ell+1-t)!}    
(2\ell+1) & \text{for all} \;\; t \in [\ell+1;2\ell+1].
  \end{cases}
\end{equation}


Recall that for odd $N$, by definition,
$N!!$ is the product of all odd numbers 
from 1 to $N$.

Note that factor $2\ell+1$  is canceled 
in the first and last formulas 
for $t=0$ and $t = 2\ell+1$, respectively. 
Yet, for $t=0$  we had to replace 
$2^{3t-1} = 1/2$  by $c_0 = 1$, 
while $c_{2\ell + 1} = 2^{4\ell}$. 


\begin{theorem}
Polynomials $a_t(\ell) ,b_t(\ell)$, and $c_t(\ell)$  
have strictly positive integer coefficients and satisfy 
identity (\ref{eq-ABC-even}). 
\end{theorem}

\proof  
Polynomials $A,B$, and $C$ are given by explicit formulas, 
so identities $k_{2\ell} + 1 = kA^2+1 = (k+1)B^2 = C$
are directly verifiable.\footnote{Note that the corresponding, pretty complicated,  
binomial identities are missing in the famous collection by Gould \cite{Gou72}.}

By the above definitions, all coefficients are strictly positive, and  $b_t$ are integer. 
Then, $a_t$ are integer, by Observation \ref{obs-integer-even}, and $c_t$  are integer, by (\ref{eq-ABC-even}).
\qed 

\medskip 

Let us illustrate this theorem by several examples. The identity 
$k A_\ell^2 + 1 = (k+1) B_\ell^2$ 
means that the coefficients  \(c_t\) of \(k^t\) in these two polynomials (of the same degree \(2\ell + 1\)) are equal for all \(t\).
For \(t \leq \ell\), we have:

\begin{align*}
t=0: \quad & 1 = b_0^2; \\[6pt]
t=1: \quad & a_0^2 = 2 b_0 b_1 + b_0^2; \\[6pt]
t=2: \quad & 2 a_0 a_1 = 2 b_0 b_2 + 2 b_0 b_1 + b_1^2; \\[6pt]
t=3: \quad & 2 a_0 a_2 + a_1^2 = (2 b_0 b_2 + 2 b_0 b_3) + (2 b_1 b_2 + b_1^2); \\[6pt]
t=4: \quad & 2 a_0 a_3 + 2 a_1 a_2 = (2 b_0 b_3 + 2 b_0 b_4) + (2 b_1 b_2 + 2 b_1 b_3) + b_2^2; \\[6pt]
t=5: \quad & 2 a_0 a_4 + 2 a_1 a_3 + a_2^2 = (2 b_0 b_4 + 2 b_0 b_5) + (2 b_1 b_3 + 2 b_1 b_4) + (2 b_2 b_3 + b_2^2); \dots
\end{align*}

For general \(t\), writing \(t = 2u\) or \(t = 2u + 1\), we have:

\begin{align*}
t = 2u: \quad & 2 a_0 a_{t-1} + 2 a_1 a_{t-2} + \dots + 2 a_{u-1} a_u \\
&= (2 b_0 b_{t-1} + 2 b_0 b_t) + (2 b_1 b_{t-2} + 2 b_1 b_{t-1}) + \dots + (2 b_{u-1} b_u + b_u^2); \\[6pt]
t = 2u + 1: \quad & 2 a_0 a_{t-1} + 2 a_1 a_{t-2} + \dots + a_u^2 \\
&= (2 b_0 b_t + 2 b_0 b_{t-1}) + (2 b_1 b_{t-1} + 2 b_1 b_{t-2}) + \dots + (2 b_u b_{u+1} + b_u^2).
\end{align*}

\smallskip 
\noindent 
These identities are defined by the following rule:  
the sum of two subscripts in a pair is \(t - 1\) for \(A_\ell\),  
and it is \(t - 1\) or \(t\) for \(B_\ell\).

\medskip

Let us verify these identities and compute \(c_t\) for \(t \leq 3\).  
To do so, “simplify” the expressions for \(a_t\) and \(b_t\)  
by cancelling factorials and powers of~2:

\[
\begin{aligned}
a_0 &= 2\ell + 1, \quad
a_1 = \frac{2}{3} \ell (\ell + 1)(2\ell + 1), \quad
a_2 = \frac{2}{15} (\ell - 1)\ell(\ell + 1)(\ell + 2)(2\ell + 1); \\[6pt]
b_0 &= 1, \quad
b_1 = 2\ell(\ell + 1), \quad
b_2 = \frac{2}{3} (\ell - 1)\ell(\ell + 1)(\ell + 2), \\[6pt]
b_3 &= \frac{4}{45} (\ell - 2)(\ell - 1)\ell(\ell + 1) \cdot (\ell + 2)(\ell + 3).
\end{aligned}
\]


\medskip

\noindent
For $t = 0$, the identity $b_0^2 = 1$ implies that $b_0 = c_0 = 1$ too.

\medskip

\noindent
For $t = 1$, the identity  $a_0^2 = 2b_0b_1 + b_0^2$ holds, since $(2\ell + 1)^2 = 4\ell(\ell + 1) + 1 = c_1$.

\medskip

\noindent
For $t = 2$, the identity $2a_0a_1 = b_1^2 + 2b_0b_1 + 2b_0b_2$ means that
\[
\frac{4}{3}(2\ell + 1)^2\ell(\ell + 1) =
4\ell^2(\ell + 1)^2 + 4\ell(\ell + 1) +
\frac{4}{3}(\ell - 1)\ell(\ell + 1)(\ell + 2).
\]


\noindent
Factoring out $\frac{4}{3}\ell(\ell + 1)$, we obtain:
$(2\ell + 1)^2 =
3\ell(\ell + 1) + 3 + (\ell - 1)(\ell + 2)$.
This can be rewritten as 
$4\ell^2 + 4\ell + 1 = 3\ell^2 + 3\ell + 3 + \ell^2 + \ell - 2$,
which holds with  $c_2 = \frac{4}{3}\ell(\ell + 1)(2\ell + 1)^2$.

\medskip

\noindent
For $t = 3$, the identity $2a_0a_2 + a_1^2 = 2b_3 + 2b_2 + 2b_1b_2 + b_1^2$
turns into
\[
\frac{4}{15}(2\ell + 1)^2(\ell - 1)\ell(\ell + 1)(\ell + 2) +
\frac{4}{9}(2\ell + 1)^2\ell^2(\ell + 1)^2 =
\]

\[
\frac{8}{45}(\ell - 2)(\ell - 1)\ell(\ell + 1)(\ell + 2)(\ell + 3) +
\frac{4}{3}(\ell - 1)\ell(\ell + 1)(\ell + 2) +
\]
\[
\frac{8}{3}(\ell - 1)\ell^2(\ell + 1)^2(\ell + 2) +
4\ell^2(\ell + 1)^2.
\]

\smallskip

\noindent
Factoring out $\frac{4}{45}\ell(\ell + 1)$, we obtain:
\[
3(2\ell + 1)^2(\ell - 1)(\ell + 2) +
5(2\ell + 1)^2\ell(\ell + 1) =
\]
\[
2(\ell - 2)(\ell - 1)(\ell + 2)(\ell + 3) +
15(\ell - 1)(\ell + 2) +
30(\ell - 1)\ell(\ell + 1)(\ell + 2) +
45\ell(\ell + 1).
\]

\smallskip

\noindent
This can be rewritten as
\[
(12\ell^4 + 24\ell^3 - 9\ell^2 - 21\ell - 6) +
(20\ell^4 + 40\ell^3 + 25\ell^2 + 5\ell) =
\]
\[
(2\ell^4 + 4\ell^3 - 14\ell^2 - 16\ell + 24) +
(30\ell^4 + 60\ell^3 - 30\ell^2 - 60\ell) +
\]
\[
(45\ell^2 + 45\ell) + (15\ell^2 + 15\ell - 30),
\]
which holds with
\[
c_3 = \frac{8}{45} \ell(\ell + 1)(2\ell - 1)(2\ell + 1)^2(2\ell + 3).
\]

\medskip 

Consider now the case $\ell < t \leq 2\ell+1$ and 
set $\tau = 2\ell+1 - t$. Then, obviously, $0 \leq \tau \leq \ell$.

\medskip 

If $t = 2\ell+1$ (and $\tau = 0$) then 
$a_\ell^2 = b_\ell^2$, that is, 
$a_\ell = b_\ell$, since both are nonnegative.

\smallskip 

If $t = 2\ell$ (and $\tau = 1$) then 
$2a_\ell a_{\ell-1} = 2b_\ell b_{\ell-1} + b_\ell^2$.

\smallskip 

If $t = 2\ell-1$ (and $\tau = 2$) then 
$2a_\ell a_{\ell-2} + a_{\ell-1}^2 = 
2b_\ell b_{\ell-1} + 2b_\ell b_{\ell-2} + b_{\ell-1}^2$.

\smallskip 

If $t = 2\ell-2$ (and $\tau = 3$) then 
$2a_\ell a_{\ell-3} + 2a_{\ell-1} a_{\ell-2} = 
2b_\ell b_{\ell-2} + 2b_\ell b_{\ell-3} + b_{\ell-1}^2$.

\smallskip 

If $t = 2\ell-3$ (and $\tau = 4$) then 
\begin{align*}
2a_\ell a_{\ell-4} + 2a_{\ell-1} a_{\ell-3} + a_{\ell-2}^2 
&= 2b_\ell b_{\ell-3} + 2b_\ell b_{\ell-4} 
+ 2b_{\ell-1} b_{\ell-2} + b_{\ell-2}^2.
\end{align*}

\smallskip 

If $\tau = 2u$ is even and, respectively, $t = 2(\ell-u)+1$ is odd then 
\begin{align*}
&2a_\ell a_{\ell-2u} + 2a_{\ell-1} a_{\ell-2u+1} + \dots + 
2a_{\ell-u+1} a_{\ell-u-1} + a_{\ell-u}^2 \\
&= 2b_\ell b_{\ell-2u+1} + 2b_\ell b_{\ell-2u} + 
2b_{\ell-1} b_{\ell-2u+2} + 2b_{\ell-1} b_{\ell-2u+1} + \dots 
+ 2b_{\ell-u+1} b_{\ell-u} + b_{\ell-u}^2.
\end{align*}

\smallskip 

If $\tau = 2u+1$ is odd and, respectively, $t = 2(\ell-u)$ is even then 
\begin{align*}
&2a_\ell a_{\ell-2u-1} + 2a_{\ell-1} a_{\ell-2u} + 
2a_{\ell-2} a_{\ell-2u+1} + \dots + 
2a_{\ell-u} a_{\ell-u-1} \\
&= 2b_\ell b_{\ell-2u+1} + 2b_\ell b_{\ell-2u} + 
2b_{\ell-1} b_{\ell-2u+2} + 2b_{\ell-1} b_{\ell-2u+1} + \dots 
+ 2b_{\ell-u+1} b_{\ell-u} + b_{\ell-u}^2.
\end{align*}

\medskip 

These identities are defined by a similar rule:  
the sum of two subscripts in a pair  
is $t - 1$ for $A_\ell$ and it is $t - 1$ or $t$ for $B_\ell$.  
Yet, all terms with a subscript greater than $\ell$  
are waved.  

\medskip

Let us verify the above identities and compute $c_t$  
for several small values of $\tau = 2\ell + 1 - t$.

\medskip 

For $\tau = 0$ we have  
$a_\ell = b_\ell = 2^{2\ell}$, and  
$c_{2\ell+1} = a^2_\ell = b^2_\ell = 2^{4\ell}$.

\medskip

For $\tau = 1$ we further have  
$\;\; a_{\ell-1} = 2^{2(\ell-1)}(2\ell+1)$, \quad
$b_{\ell-1} = 2^{2(\ell-1)}(2\ell-1)$,  
\[
c_{2\ell} = 2a_\ell a_{\ell-1} = 
2b_\ell b_{\ell-1} + b_\ell^2 = 
2^{4\ell-1}(2\ell+1).
\]

\medskip

For $\tau = 2$ we have  
$\;\; a_{\ell-2} = 2^{2\ell-4}(\ell-1)(2\ell+1)$, \quad
$b_{\ell-2} = 2^{2\ell-4}(2\ell-3)(\ell-1)$,  
\[
c_{2\ell-1} = 2a_\ell a_{\ell-2} + a_{\ell-1}^2 =
(2b_\ell b_{\ell-1} + 2b_\ell b_{\ell-2}) + b_{\ell-1}^2 =
2^{4\ell-4}(2\ell+1)(4\ell-1).
\]

\medskip

For $\tau = 3$ we have  
\[
\begin{aligned}
a_{\ell-3} &= \frac{2^{2\ell-6}}{3}(2\ell+1)(\ell-2)(2\ell-3), \\
b_{\ell-3} &= \frac{2^{2\ell-8}}{3}(2\ell-5)(\ell-2)(2\ell-3),
\end{aligned}
\]
\[
\begin{aligned}
c_{2\ell-2} &= 2a_\ell a_{\ell-3} + 2a_{\ell-1} a_{\ell-2} 
= (2b_\ell b_{\ell-3} + 2b_\ell b_{\ell-2}) +
(2b_{\ell-1} b_{\ell-2} + b_{\ell-1}^2) \\
&= \frac{2^{4\ell-5}}{3}(2\ell+1)(2\ell-1)(4\ell-3).
\end{aligned}
\]

\medskip

For $\tau = 4$ we have  
\[
\begin{aligned}
a_{\ell-4} &= \frac{2^{2\ell-9}}{3}(2\ell+1)(2\ell-5)(\ell-2)(\ell-3), \\
b_{\ell-4} &= 2^{2\ell-9}(2\ell-7)(\ell-3)(2\ell-5)(\ell-2),
\end{aligned}
\]
\[
\begin{aligned}
c_{2\ell-3} &= 2a_\ell a_{\ell-4} + 2a_{\ell-1} a_{\ell-3} + a_{\ell-2}^2 
= (2b_\ell b_{\ell-3} + 2b_\ell b_{\ell-4}) +
(2b_{\ell-1} b_{\ell-2} + 2b_{\ell-1} b_{\ell-3}) + b_{\ell-2}^2 \\
&= \frac{2^{4\ell-8}}{3}(2\ell+1)(4\ell-5)(\ell-1)(4\ell-3).
\end{aligned}
\]

\medskip

For $\tau = 5$ we have  
\[
\begin{aligned}
a_{\ell-5} &= \frac{2^{2\ell-11}}{15}(2\ell+1)(\ell-4)(2\ell-7)(\ell-3)(2\ell-5), \\
b_{\ell-5} &= \frac{2^{2\ell-10}}{15}(2\ell+1)
(4\ell-6)(\ell-2)(4\ell-7)(2\ell-3)(4\ell-5),
\end{aligned}
\]

\[
\begin{aligned}
c_{2\ell-4} &= 2a_\ell a_{\ell-5} + 2a_{\ell-1} a_{\ell-4} + 2a_{\ell-2} a_{\ell-3} \\& 
= (2b_\ell b_{\ell-5} + 2b_\ell b_{\ell-4}) +
(2b_{\ell-1} b_{\ell-4} + 2b_{\ell-1} b_{\ell-3}) +  
(2b_{\ell-2} b_{\ell-3} + b_{\ell-2}^2) \\
&= \frac{2^{4\ell-6}}{15}(2\ell+1)(4\ell-7)(2\ell-3)(4\ell-5)(\ell-1).
\end{aligned}
\]
\section{Case of odd $i$}
Given integer $\ell \geq 0, \; i = 2\ell+1$,  
and  $t$ such that 
$0 \leq t \leq \ell$ for $a$ and 
$0 \leq t \leq \ell+1$ for $b$, 
set $b_0 = 1$, 
while all others $a_t$ and $b_t$ are defined 
by the next  
chains of identical expressions: 
\begin{equation} 
\begin{aligned}
a_t & = 2^{2t+1} \frac{(\ell+t+1)!}{(2t+1)!(\ell-t)!} = 
2^{2t+1}\binom{\ell+t+1}{\ell-t} = 
2^{2t+1}\binom{\ell+t+1}{2t+1},\\     
b_t & = 2^{2t-1} \frac{(\ell+t)!}{(2t-1)! (\ell-t+1)!}  
\frac{\ell+1}{t} =
2^{2t-1} \binom{\ell+t}{2t-1} \frac{\ell+1}{t}  = 
2^{2t-1} \binom{\ell+t}{\ell-t+1} \frac{\ell+1}{t}.
\end{aligned}
\end{equation}

Obviously, $a_t$ is integer, 
due to its binomial expression. 

\begin{observation}
\label{obs-integer-odd}
Coefficients $b_t$  are integer too.
\end{observation}

\proof 
Again, as in the even case,  we will show more, namely, that 
$\frac{b_t}{2^{2t+1}} = 
\binom{\ell+t}{\ell-t+1}  
\frac{\ell+1}{t}$  are integer. 
Note that 
$\frac{(\ell-t+2) \times \dots \times (\ell+t)}
{2 \times \dots \times (2t-1)}$  is a multiple of $t$ 
unless $\ell + 2 \equiv 1 \pmod t$, 
but in this case $\ell+1$ is a multiple of $t$.
\qed 

\medskip 

Introduce polynomial 
$C_i = \sum_{t=0}^{2(\ell+1)} c_t k^t$, where $i = 2\ell+1$, 
by the following formula: 

\begin{equation}
  c_t = 
  \begin{cases}
    \frac{2^{3t-1}}{(2t)!} 
\frac{(\ell+\frac{t}{2})!}{(\ell-\frac{t}{2}+1)!}  
\frac{(2\ell+t+1)!!}{(2\ell-t+1)!!}(\ell+1)  & \text{for even} \;\; t \in [2;\ell], \text{while} \;\; c_0 = 1; \\ 
    \frac{2^{3t}}{(2t)!} 
\frac{(\ell+\frac{t+1}{2})!}{(\ell-\frac{t-1}{2})!}  
\frac{(2\ell+t)!!}{(2\ell-t+2)!!}(\ell+1)  & 
\text{for odd} \;\; t \in [1;\ell];\\
2^{2t} 
\frac{(2\ell+t+1)!}{(2t)!(2\ell+2-t)!} (\ell+1)  &\text{for all} \;\; t \in [\ell+1;2(\ell+1)].
  \end{cases}
\end{equation}

Note that in the first and last formulas, 
for $t=0$ and $t = 2(\ell+1)$, 
factor $\ell+1$ is cancelled and  
we obtain $c_{2(\ell+1)} = 2^{4\ell+2}$, which is correct, 
yet, for $t=0$ we had to replace
$c_t = 2^{3t-1} = \frac{1}{2}$ by $c_0 = 1$. 

\begin{theorem}
Polynomials $a_t(\ell) ,b_t(\ell)$, and $c_t(\ell)$  
have strictly positive integer coefficients  and satisfy 
identity (\ref{eq-ABC-odd}).    
\end{theorem}

\proof  
We can just repeat the proof of the previous theorem  
(and the footnote). 
\newline 
Polynomials $A,B$, and $C$ are given by explicit formulas, 
so identities  
\newline 
$k_{2\ell+1} + 1 = k(k+1)A^2+1 = B^2 = C$
are directly verifiable. 
\qed 

\medskip 

Again we illustrate this theorem by several examples.  
Identity \( k(k+1)A_\ell^2+1 = B_{\ell+1}^2 \) holds if and only if 
these two polynomials have the same coefficients \( c_t \) for each term \( k^t \), 
for all \( t = 0,1,\dots,\ell \).  
Thus, we have to check the following binomial identities:
\begin{flalign*}
& t=0: \quad 1 = b_0^2 && \text{(even \( t \))} & \\
& t=1: \quad a_0^2 = 2b_0b_1 && \text{(odd \( t \))} & \\
& t=2: \quad 2a_0 a_1 + a_0^2 = 2b_0b_2 + b_1^2 && \text{(even \( t \))} & \\
& t=3: \quad 2a_0a_2 + 2a_0a_1 + a_1^2 = 2b_0b_3 + 2b_1b_2 && \text{(odd \( t \))} & \\
& t=4: \quad (2a_0a_3 + 2a_0a_2) + (2a_1a_2 + a_1^2) = 2b_0b_4 + 2b_1 b_3 + b_2^2 && \text{(even \( t \))} & \\
& t=5: \quad (2a_0a_4 + 2a_0a_3) + (2a_1a_3 + 2a_1a_2) + a_2^2 = 2b_0b_5 + 2b_1b_4 + 2b_2b_3 && \text{(odd \( t \))} & \\
& t=2u: \quad (2a_0 a_{t-1} + 2a_0 a_{t-2}) + (2a_1a_{t-2} + 2a_1a_{t-3}) + \dots + (2a_u a_{u-1} + a_u^2) = {} &&& \\
& \quad\quad 2b_0b_t + 2b_1 b_{t-1} + \dots + b_u^2 && \text{(even \( t \))} & \\
& t=2u+1: \quad (2a_0 a_{t-1} + 2a_0 a_{t-2}) + (2a_1a_{t-2} + 2a_1a_{t-3}) + \dots + a_u^2 = {} &&& \\
& \quad\quad 2b_0 b_t + 2b_1 b_{t-1} + \dots + 2b_u b_{u+1} && \text{(odd \( t \))} &
\end{flalign*}

\medskip

These identities are defined by the rule: 
the sum of two subscripts in a pair 
is \( t \) for \( B_{\ell+1} \) and it is \( t-2 \) or \( t-1 \) for \( A_\ell \).  
Let us verify these identities and compute polynomial 
\( c_t = c_t(\ell) \) for \( t \leq 2 \).  
To do so, let us ``simplify'' \( a_t \) and \( b_t \) by canceling factorials and 
powers of~2: 

\medskip

\[
\begin{aligned}
a_0 &= 2(\ell+1), \\
a_1 &= \frac{4}{3}\,\ell(\ell+1)(\ell+2)(2\ell+1), \\
a_2 &= \frac{4}{15}\,(\ell-1)\ell(\ell+1)(\ell+2)(2\ell+1), \\
a_3 &= \frac{8}{315}(\ell-2)(\ell-1)\ell(\ell+1)(\ell+2)(\ell+3)(\ell+4)(2\ell+1) \\
\\[0.1em]
b_0 &= 1, \\
b_1 &= 2(\ell+1)^2, \\
b_2 &= \frac{2}{3}\,\ell(\ell+1)^2(\ell+2), \\
b_3 &= \frac{4}{45}(\ell-2)(\ell-1)\ell(\ell+1)(\ell+2)(\ell+3)
\end{aligned}
\]

\medskip

For \( t=0 \), the identity is reduced to \( 1 = b_0^2 \), which holds (even \( t \)).

\medskip

For $t=1$  we have  
$a_0^2 = 2b_0b_1 = [2(\ell+1)]^2 = c_1.$  

\medskip 

For $t=2$ identity 
$2a_0 a_1 +  a_0^2 = 2b_0b_2 + b_1^2$ 
can be rewritten as 

$\frac{16}{3}\ell(\ell+1)^2(\ell+2) + 4(\ell+1)^2 = 
\frac{4}{3}\ell(\ell+1)^2(\ell+2) + 4(\ell+1)^4$. 

\smallskip 

Factoring $\frac{4}{3}(\ell+1)^2$ out we obtain  
$4 \ell(\ell+2) + 3 = \ell(\ell+2) + 3(\ell+1)^2$, or equivalently 

\smallskip

$3 \ell(\ell+2)  + 3 = 3(\ell+1)^2, \;$ or  
$\ell(\ell+2)  + 1 = (\ell+1)^2,$ which holds too, 
with 

$c_2 = \frac{4}{3}(\ell+1)^2(2\ell+1)(2\ell+3)$.

\bigskip 


\bigskip 

Assume now that $\ell + 1 < t \leq 2(\ell + 1)$ and set $\tau = 2\ell + 2 - t$, 
so that $0 \leq \tau \leq \ell$. Consider successively the values $\tau = 0, 1, \dots$. 
The identities for $c_t$ are defined by a similar rule: 
the sum of the two subscripts in a pair is $t - 1$ or $t - 2$ for $A_\ell$, and 
it is exactly $t$ for $B_{\ell + 1}$. 
Yet, all terms with subscripts greater than $\ell$ for $A$ and greater than $\ell + 1$ for $B$ are omitted.

Let us verify these identities and compute the polynomial $c_t = c_t(\ell)$ for $\ell < t \leq 2\ell + 2$. 
To do so, ``simplify'' $a_t$ and $b_t$ for $\ell + 1 \geq t \geq \ell - 5$ by canceling factorials and powers of 2:

\medskip

\[
\begin{aligned}
a_\ell &= 2^{2\ell+1}, \\
a_{\ell-1} &= 2^{2\ell} \ell, \\
a_{\ell-2} &= 2^{2\ell-3} (\ell-1)(2\ell-1), \\
a_{\ell-3} &= \frac{2^{2\ell-4}}{3} (\ell-2)(2\ell-3)(\ell-1), \\
a_{\ell-4} &= \frac{2^{2\ell-8}}{3} (\ell-3)(2\ell-5)(\ell-2)(2\ell-3), \\
a_{\ell-5} &= \frac{2^{2\ell-9}}{45} (\ell-4)(2\ell-7)(\ell-3)(2\ell-5)(\ell-2).
\end{aligned}
\]

\medskip

\[
\begin{aligned}
b_{\ell+1} &= 2^{2\ell+1}, \\
b_\ell &= 2^{2\ell}(\ell+1), \\
b_{\ell-1} &= 2^{2\ell-3} (\ell+1)(2\ell-1), \\
b_{\ell-2} &= \frac{2^{2\ell-4}}{3} (\ell+1)(2\ell-3)(\ell-1), \\
b_{\ell-3} &= \frac{2^{2\ell-8}}{3} (\ell+1)(2\ell-5)(\ell-2)(2\ell-3), \\
b_{\ell-4} &= \frac{2^{2\ell-9}}{3} (\ell+1)(2\ell-7)(\ell-3)\ell(2\ell-5)(\ell-2), \\
b_{\ell-5} &= \frac{2^{2\ell-12}}{45} (\ell+1)(2\ell-9)(\ell-4)(2\ell-7)(\ell-3)\ell(2\ell-5).
\end{aligned}
\]

\medskip

For $t = 2\ell + 2$ (and $\tau = 0$), we must have $\;\; a_\ell^2 = b_{\ell+1}^2,$
that is, $a_\ell = b_{\ell+1}$ since both are nonnegative.  
Since $b_{\ell+1} = a_\ell = 2^{2\ell+1}$, this equality holds and $\; c_{2\ell+2} = 2^{4\ell+2}$.

\smallskip

For $t = 2\ell + 1$ (and $\tau = 1$), the equality  $\; 2 b_{\ell+1} b_\ell = 2 a_\ell a_{\ell-1} + a_\ell^2$ 
holds with 

$c_{2\ell+1} = 2^{4\ell+2} (\ell+1).$

\smallskip

For $t = 2\ell$ (and $\tau = 2$), the equality
$\; 2 b_{\ell+1} b_{\ell-1} + b_\ell^2 = (2 a_\ell a_{\ell-1} + 2 a_\ell a_{\ell-2}) + a_{\ell-1}^2$ 
\newline 
holds with
$\; c_{2\ell} = 2^{4\ell - 1} (\ell+1)(4\ell + 1)$.

\smallskip

For $t = 2\ell - 1$ (and $\tau = 3$), the equality
\[
2 b_{\ell+1} b_{\ell-2} + 2 b_\ell b_{\ell-1} = (2 a_\ell a_{\ell-2} + 2 a_\ell a_{\ell-3}) + (2 a_{\ell-1} a_{\ell-2} + a_{\ell-1}^2)
\]
holds with
\[
c_{2\ell - 1} = \frac{2^{4\ell - 1}}{3} (\ell+1)(4\ell - 1) \ell.
\]

\smallskip

For $t = 2\ell - 2$ (and $\tau = 4$), the equality
\[
2 b_{\ell+1} b_{\ell-3} + 2 b_\ell b_{\ell-2} + b_{\ell-1}^2 = (2 a_\ell a_{\ell-4} + 2 a_\ell a_{\ell-3}) + (2 a_{\ell-1} a_{\ell-3} + 2 a_{\ell-1} a_{\ell-2})
\]
holds with
\[
c_{2\ell - 2} = \frac{2^{4\ell - 6}}{3} (\ell+1) (4\ell - 3)(2\ell - 1)(4\ell - 1).
\]

\smallskip

For $t = 2\ell - 3$ (and $\tau = 5$), the equality
\[
\begin{aligned}
2 b_{\ell+1} b_{\ell-4} + 2 b_\ell b_{\ell-3} + 2 b_{\ell-1} b_{\ell-2} = {} & \\
(2 a_\ell a_{\ell-5} + 2 a_\ell a_{\ell-4}) + {} & (2 a_{\ell-1} a_{\ell-4} + 2 a_{\ell-1} a_{\ell-3}) + 
(2 a_{\ell-2} a_{\ell-3} + a_{\ell-2}^2)
\end{aligned}
\]
holds with
\[
c_{2\ell - 3} = \frac{2^{4\ell - 9}}{15} (\ell+1) (4\ell - 5)(\ell - 1)(4\ell - 3)(2\ell - 1).
\]

\medskip 

\section{Identities for $i \leq 8$}
The polynomials for \(k_i\), \(k_i + 1\), and the corresponding square quadratic numbers
\newline 
$\Pi_i = P(k, k_i), \quad 1 \leq i \leq 8$,
are as follows:

\medskip

\[
\begin{aligned}
k_1 + 1 &= 4k(k+1) + 1 = (2k+1)^2, \\
\Pi_1 &= \bigl[2k(k+1)(2k+1)\bigr]^2;
\end{aligned}
\]

\medskip

\[
\begin{aligned}
k_2 + 1 &= k(4k+3)^2 + 1 = (k+1)(4k+1)^2, \\
\Pi_2 &= \bigl[k(k+1)(4k+1)(4k+3)\bigr]^2;
\end{aligned}
\]

\medskip

\[
\begin{aligned}
k_3 + 1 &= 16k(k+1)(2k+1)^2 + 1 = (8k^2 + 8k + 1)^2, \\
\Pi_3 &= \bigl[4k(k+1)(2k+1)(8k^2 + 8k + 1)\bigr]^2;
\end{aligned}
\]

\medskip

\[
\begin{aligned}
k_4 + 1 &= k \bigl((4k)^2 + 5(4k) + 5\bigr)^2 + 1 = (k+1) \bigl((4k)^2 + 3(4k) + 1\bigr)^2, \\
\Pi_4 &= \bigl[k(k+1)((4k)^2 + 3(4k) + 1)((4k)^2 + 5(4k) + 5)\bigr]^2;
\end{aligned}
\]

\medskip

\[
\begin{aligned}
k_5 + 1 &= 4k(k+1)(4k+1)^2(4k+3)^2 + 1 = (2k+1)^2 (16k^2 + 16k + 1)^2, \\
\Pi_5 &= \bigl[2k(k+1)(2k+1)(4k+1)(4k+3)(16k^2 + 16k + 1)\bigr]^2;
\end{aligned}
\]

\medskip

\[
\begin{aligned}
k_6 + 1 &= k \bigl((4k)^3 + 7(4k)^2 + 14(4k) + 7\bigr)^2 + 1 = (k+1) \bigl((4k)^3 + 5(4k)^2 + 6(4k) + 1\bigr)^2, \\
\Pi_6 &= \bigl[k(k+1)((4k)^3 + 5(4k)^2 + 6(4k) + 1)((4k)^3 + 7(4k)^2 + 14(4k) + 7)\bigr]^2;
\end{aligned}
\]

\medskip

\[
\begin{aligned}
k_7 + 1 &= 64 k(k+1)(16k^3 + 24k^2 + 10k + 1)^2 + 1 \\
&= 64 k(k+1)(2k+1)^2 (8k^2 + 8k + 1)^2 + 1 \\
&= \bigl(128 k^4 + 256 k^3 + 160 k^2 + 32 k + 1\bigr)^2, \\
\Pi_7 &= \bigl[k(k+1)(2k+1)(8k^2 + 8k + 1)(128 k^4 + 256 k^3 + 160 k^2 + 32 k + 1)\bigr]^2;
\end{aligned}
\]

\medskip

\[
\begin{aligned}
k_8 + 1 &= k \bigl((4k)^4 + 9(4k)^3 + 27(4k)^2 + 30(4k) + 9\bigr)^2 + 1 \\
&= (k+1) \bigl((4k)^4 + 7(4k)^3 + 15(4k)^2 + 10(4k) + 1\bigr)^2, \\
\Pi_8 &= \Bigl[k(k+1)((4k)^4 + 7(4k)^3 + 15(4k)^2 + 10(4k) + 1) \\
&\quad \times \bigl((4k)^4 + 9(4k)^3 + 27(4k)^2 + 30(4k) + 9\bigr)\Bigr]^2.
\end{aligned}
\]

\medskip 

For example, (\ref{eq-odd-new}), in case $i=3$, is reduced to: 

\smallskip 

$k_3 + 1 = 
k(k+1)(a_1 k + a_0)^2 + 1 = (b_2k^2 + b_1k + b_0)^2$. 

\smallskip 

Equating coefficients for $k^t$  for  $t=0,1,2,3,4$, we obtain: 

\smallskip    

$b_0=1, \; a_0^2 = 2b_1 b_0$,  
$2a_1 a_0 + a^2_0 = 2b_2 b_0 + b^2_1, \;$ 
$2a_1 a_0 + a^2_1 = 2b_2 b_1, \;$ 
$a^2_1 = b^2_2.$

\smallskip  

Hence, $b_0 = 1, b_1 = a^2_0/2$, and $b_2 = \pm a_1$. 
In case $b_2 = -a_1$  there are no integer solutions,  
while setting $b_2 = a_1$, we obtain equations 
$2a_1(a_0 - 1) = a^2_0/4 - a^2_0$ and 
$a^2_1 = a^2_0 - 2a_0$, 
the last one, assuming that $a_1 \neq 0.$ 
Then, for $a_0$ we obtain the cubic equation  
$a^3_0 - 8a^2_0 + 20 a_0 - 16 = 0$ 
with roots: $a_0=2$, of multiplicity two, and $a_0=4$.
The first one results in 
$a_1=b_2=0$ and $a_0 = b_1=2$, 
in contradiction to our assumption: $a_1 \neq 0$.
(Note that it 
results in  
$4k(k+1) + 1 = (2k+1)^2 = k_1 + 1$.)  
The second one results in 
$a_0=4, a_1=b_1=b_2=8$, that is, 
$8k(k+1)(2k+1)^2+1 = (8k^2 +8k+1)^2=k_3 +1$. 

\medskip 

For $t=1$ we have $k_1 - k = k(4k+3)$. 
Furthermore, $k_1 = 4k(k+1), \; k+j+1 = (2k+1)^2$, 
Hence, $\Pi_1(k) = k(k+1)k_1(k_1+1) = [2k(k+1)(2k+1)]^2$.  
 is a perfect square for all integer $k$ with  
$|\sqrt{\Pi(k)}| = 2k(k+1)(2k+1)$. 

For all positive integer $k$ and $i \leq 8$  we obtain similar identities:

\medskip 

For \( k_2 - k = k[(4k)^2 + 3(4k) + 8] = 8k(k+1)(2k+1) \) we have:

\smallskip 

\[
k_2 = k(4k+3)^2 \quad \text{and} \quad k_2 + 1 = (k+1)(4k+1)^2.
\]
Hence, the quadratic number \(\Pi_2(k)\) is a perfect square for all \(k\), with
\[
\sqrt{\Pi_2} = k(k+1)(4k+1)(4k+3).
\]

\bigskip

For \(k_3 - k = k(4k+3)(16k^2 + 20k + 5) \) we have:

\[
k_3 = 16k(k+1)(2k+1)^2 = 8 \sqrt{\Pi_1}, \quad k_3 + 1 = (8k^2 + 8k + 1)^2.
\]

Hence,
\[
\sqrt{\Pi_3} = 4k(k+1)(2k+1)(8k^2 + 8k + 1).
\]

\bigskip

For \( k_4 - k = 8k(k+1)(2k+1)(16k^2 + 16k + 3) \) we have:

\[
k_4 = k \bigl[(4k)^2 + 5(4k) + 5\bigr]^2, \quad
k_4 + 1 = (k+1) \bigl[(4k)^2 + 3(4k) + 1\bigr]^2.
\]

Hence, \(\Pi = \Pi(k, k_4)\) is a perfect square for all \(k\), with
\[
\sqrt{\Pi_4} = k(k+1)(16k^2 + 12k + 1)(16k^2 + 20k + 5).
\]

\medskip

For \( k_5 -k = k(16k^2 + 20k + 5)(64k^3 + 112k^2 + 56k + 7) \) we have

\[
k_5 = 4k(k+1)(4k+1)^2 (4k+3)^2, \quad
k_5 + 1 = (2k+1)^2 (16k^2 + 16k + 1)^2.
\]

Hence, \(\Pi_5(k)\) is a perfect square for all integer \(k\), with
\[
\sqrt{\Pi_5} = 2k(k+1)(2k+1)(4k+1)(4k+3)(16k^2 + 16k + 1).
\]

\medskip

For
\[
k_6 - k = k \left( \bigl[(4k)^3 + 7(4k)^2 + 14(4k) + 7\bigr]^2 - 1 \right) = 16k(k+1)(2k+1)(4k+1)(4k+3)(8k^2 + 8k + 1),
\]
we have
\[
k_6 = k \bigl[(4k)^3 + 7(4k)^2 + 14(4k) + 7\bigr]^2,
\quad
k_6 + 1 = (k+1) \bigl[(4k)^3 + 5(4k)^2 + 6(4k) + 1\bigr]^2,
\]
and
\[
\sqrt{\Pi_6} = \sqrt{\Pi(k,k_6} = k(k+1) \bigl[(4k)^3 + 5(4k)^2 + 6(4k) + 1\bigr] \bigl[(4k)^3 + 7(4k)^2 + 14(4k) + 7\bigr].
\]

\medskip

For
\[
k_7 -  k = \bigl[64(k+1)(16k^3 + 24k^2 + 10k + 1)^2 - 1\bigr],
\]
\[
k_7 = 64 k (k+1)(16k^3 + 24k^2 + 10k + 1)^2 = 64 k (k+1)(2k+1)^2 (8k^2 + 8k + 1)^2,
\]
\[
k_7 + 1 = (128 k^4 + 256 k^3 + 160 k^2 + 32 k + 1)^2,
\]
and
\[
\sqrt{\Pi_7} = \sqrt{\Pi(k,k_7)} = k(k+1)(2k+1)(8k^2 + 8k + 1)(128 k^4 + 256 k^3 + 160 k^2 + 32 k + 1).
\]

\medskip

Finally, for
\[
k_8 = k \left[ \left( (4k)^4 + 9(4k)^3 + 27(4k)^2 + 30(4k) + 9 \right)^2 - 1 \right]
\]
and
\[
\sqrt{\Pi_8} = \sqrt{\Pi(k,k_8)} = k(k+1) \bigl[(4k)^4 + 7(4k)^3 + 15(4k)^2 + 10(4k) + 1\bigr] \times 
\]
\[
\quad \times \bigl[(4k)^4 + 9(4k)^3 + 27(4k)^2 + 30(4k) + 9\bigr].
\]


\section {Generalizations and extensions} 
We aim to characterize positive integer $s$-vectors 
\[
\mathbf{k} = (k_1, \dots, k_s)
\]
such that one of the following holds:

\begin{itemize}
  \item[(i)] The product $\Pi_s = \prod_{r=1}^s k_r(k_r + 1)$ is a perfect square; or
  \item[(ii)] $2\Pi_s$ is a perfect square.
\end{itemize}

\smallskip

For $s = 1$, this problem was resolved by Euler in 1778. 
It is easy to see that no solution of type (i) exists in this case. 
Indeed, the product $k(k + 1)$ cannot be a perfect square, 
because $k$ and $k + 1$ are consecutive integers and hence coprime, so they cannot both be perfect squares simultaneously.

\smallskip

For $s = 2$, a large family of solutions of type (i) was constructed in the present section. 

\smallskip

In contrast to the case $s = 1$, 
both types of solutions exist when $s = 2$. For example, the identity
\[
2 \cdot (3 \cdot 4) \cdot (24 \cdot 25) = (2^3 \cdot 3 \cdot 5)^2 = 120^2
\]
shows that $(3, 24)$ is a solution of type (ii).

\smallskip

Now, consider $s \geq 2$ and an $s$-vector 
$\mathbf{k} = (k_1, \dots, k_s)$
which is a solution of either type (i) or (ii) and contains an entry equal to 2 or 3. Replace one such entry: change $2 \mapsto 3$ or $3 \mapsto 2$. It is easy to verify that the resulting vector $\mathbf{k'}$ will yield a solution of the ``opposite'' type: (i) and (ii) are interchanged.
For instance, starting from the previous type (ii) solution
$3 \cdot 4)(24 \cdot 25) = \frac{120^2}{2}$,
we replace $(3 \cdot 4$ with $2 \cdot 3$ to obtain the product
$(2 \cdot 3)(24 \cdot 25) = (3 \cdot 4 \cdot 5)^2 = 60^2$,
which is a solution of type (i).

\smallskip 

Each Euler solution of type (ii) for $s = 1$,
\[
2k(k+1) = d^2,
\]
generates a solution of type (ii) with $s = 2$ and $j = 1$ as follows:
\[
(2k(2k+1)) \cdot ((2k+1)(2k+2)) = 2[d(2k+1)]^2.
\]

\smallskip

There also exist solutions of type (ii) with $s = 2$ 
that do not arise from the above two families. 
For example:
\[
\begin{aligned}
2 \cdot (3 \cdot 4)(24 \cdot 25) &= (2^3 \cdot 3 \cdot 5)^2 = 120^2, \\
2 \cdot (4 \cdot 5)(9 \cdot 10)  &= (2^2 \cdot 3 \cdot 5)^2 = 60^2, \\
2 \cdot (5 \cdot 6)(15 \cdot 16) &= (2^3 \cdot 3 \cdot 5)^2 = 120^2, \\
2 \cdot (6 \cdot 7)(27 \cdot 28) &= (2^2 \cdot 3^2 \cdot 7)^2 = 252^2, \\
2 \cdot (7 \cdot 8)(63 \cdot 64) &= (2^5 \cdot 3 \cdot 7)^2 = 777^2, \\
2 \cdot (9 \cdot 10)(80 \cdot 81) &= (2^3 \cdot 3^3 \cdot 5)^2 = 1080^2, \\
2 \cdot (11 \cdot 12)(32 \cdot 33) &= (2^4 \cdot 3 \cdot 11)^2 = 528^2, \\
2 \cdot (15 \cdot 16)(120 \cdot 121) &= (2^4 \cdot 3 \cdot 5 \cdot 11)^2 = 2640^2.
\end{aligned}
\]

\smallskip

Characterizing all solutions of type (i) or type (ii) remains an open problem for $s \geq 2$.


\section{On pairs of intervals of successive integers with equal sum of squares}
The above results allow us to solve the problem 
if the lengths of these two intervals differ by 1. 
We reduce this problem to equation 
\begin{equation}
\label{eq1}   
(n-m)^2 + (n-m+1)^2 + \dots + (n+k)^2 = 
(n+k+j+1)^2 + (n+k+j+2)^2 + \dots + (n+2k+j)^2
\end{equation}
in integer $j, m, k$ and $n$ such that  $j \geq 0, m \geq 0$ and $k > 0$. 
We rewrite it as a quadratic equation on $n$ and obtain an explicit formula   
for its discriminant  $D = D(j,m,k)$. 
In case $m=0$  we obtain $D = \Pi(k,j)$. 
Hence, $D$ is a perfect square whenever $\Pi(k,j)$ is.   
Moreover, in this case both roots are integer. 
We also obtain some partial results for $j=0$ and 
mention shortly the known cases: $m=j=0$ and $k=1$.

\subsection{Four parameters describing two intervals of successive integers}
Identities involving sums of squares of integers, usually consecutive, 
have been studied extensively for a long time; see for example, \cite{AM11, God20}. 
Here we are looking for pairs of intervals of successive integers with equal sums of squares. 

\begin{observation}
Without loss of generality 
the problem is reduced to equation (\ref{eq1})  
in integers $j,m,k$ and $n$ 
such that $j \geq 0, m \geq 0$, and $k > 0$. 
\end{observation} 

\proof 
Equation \eqref{eq1} corresponds to the case 
of two disjoint intervals of lengths 
$k+m+1$ and $k$  with distance $j+1$ between them. 

In general, two intervals may intersect.  
In this case, we just delete 
all terms of the intersection 
from both sides of the corresponding equation. 
Clearly, we obtain an equivalent reduced equality  
that corresponds to disjoint intervals.   
Furthermore, in this new equality 
$j$ becomes equal to the size of the intersection plus 1. 
Furthermore, $j=0$ means that the two obtained disjoint integer intervals are adjacent, 
that is, at distance 1. 

Equation \eqref{eq1} assumes that 
(*) the left interval is longer than the right one. 
Their lengths differ by $m+1$, 
that is, $m\geq 0$ and if $m=0$, 
the left interval is still longer by 1.  
In general, (*) might fail. 
Then we just negate all numbers in both intervals.  
Obviously, this transformation 
respects the equality and enforces  (*). 
If both intervals are of the same length then $m=-1$.
Obviously, in this case \eqref{eq1} holds 
if and only if the considered two intervals 
are symmetric with respect to  0, that is, 
one is the negation of the other. 
We do not consider this case, since it is trivial, 
in other words, we assume that $m \geq 0$. 
\qed 

Thus, we consider two intervals, 
of  $k+m+1$  and $k$ successive integers, at distance  $j+1$, 
where $j \geq 0, m \geq 0,$ and $k > 0$ are integers. 
In contrast, $n$ may be a complex number, but we are interested 
only in integer $n$, which may be negative yet.  
We transform \eqref{eq1} to a quadratic equation 
$an^2 + bn + c = 0$  with respect to $n$  whose 
coefficients $a,b,c$ and discriminant  
$D = \sqrt{b^2 - 4ac}$ depend on $j,m,k$. 
Two solutions are given by formula 
$n = \frac{1}{2a}(-b \pm \sqrt{D})$. 
If $D < 0$, neither solution is real. 
If $D \geq 0$ both are real and they coincide if and only if $D=0$. 
Assuming $a, b \in \mathbb{Z}$, both solutions are rational if and only if $D = D(j,m,k)$ is a perfect square. Because we are interested in integer solutions, any rational candidate $n$ must have a reduced denominator dividing $2a$.

\smallskip 

We derive a formula for $D(j,m,k)$, 
which allows us to analyze \eqref{eq1} in special cases 
$m=0$ and/or $j=0$. 
The first one leads to some classical results 
on ``near-isosceles" Pythagorean triples  
and to several new results as well.  
In the last two subsections we consider 
cases $j=0$ and $0 < j \leq 18$ 
and provide some partial and computational results.

\subsection{Case $j=m=0$ and Dostor's (1879) identity} 
\label{s01}
In this case, \eqref{eq1} turns into  
$$n^2 + (n+1)^2 + \dots + (n+k)^2 = 
(n+k+1)^2 + (n+k+2)^2 + \dots + (n+2k)^2,$$ 
which is reduced to 
$$n^2 = 2k^2n + k^2(2k+1).$$ 
Interestingly, its discriminant 
$$D(0,0,k) = 4k^4 + 4k^2(2k+1) = [2k(k+1)]^2$$     
is a perfect square for all $k$, and we obtain 
the old result of Dostor \cite{Dos1879}.

\begin{theorem} 
\label{t-Dos} 
For each integer positive $k$ there exist two integer solutions 
$n' = k(2k+1)$  and  $n'' = -k$. 
The first one implies identity  

\begin{equation}
\label{eq8}
\begin{split}
[k(2k+1)]^2 + \dots + [2k(k+1)]^2 &= 
[2k(k+1) + 1]^2 + \dots + [k(2k+3)]^2 \\
& = \frac{1}{6}k(k+1)(2k+1)(12k^2+12k+1),
\end{split}
\end{equation}

\noindent 
while the second one generates 
a less exciting, but still correct, identity    
$$(-k)^2 + \dots (-1)^2 + 0^2 = 1^2 + \dots + k^2 = 
\frac{1}{6}k(k+1)(2k+1).$$ 
\qed 
\end{theorem} 

Several first solutions of the first type are 

\medskip 
\noindent 
$k=1: \; 3^2+4^2=5^2 = 25$;  

\noindent 
$k=2: \;  10^2+11^2+12^2 = 13^2+14^2 = 365$; 

\noindent 
$k=3: \;  21^2+22^2+23^2+24^2 = 25^2+26^2+27^2 = 2,030$; 
 
\noindent 
$k=4: \;36^2+37^2+38^2+39^2+40^2 = 41^2+42^2+43^2+44^2 = 7,230$; 

\noindent 
$k=5: \;55^2+56^2+57^2+58^2+59^2+60^2 = 61^2+62^2+63^2+64^2+65^2 = 19,855$; 

\noindent 
$k=6: \;  78^2+79^2+80^2+81^2+82^2+83^2+84^2 = 
85^2+86^2+87^2+88^2+89^2+90^2 = 45,955.$

\medskip 

\subsection{Transformations of equation \eqref{eq1}}
\label{s-tr}
We will rewrite \eqref{eq1}  
in several equivalent ways as follows: 
$$[(n-m)^2 + \dots + n^2] + 
[(n+1)^2 + \dots + (n+k)^2] = 
[(n+k+1+j)^2 + \dots + (n+2k+j)^2];$$
\begin{equation}
\label{eq1a}
(n-m)^2 + \dots + n^2 = 
\sum_{i=1}^k [(n+k+j+i)^2 -(n+i)^2].   
\end{equation}
Let $L$ and $R$ denote the left- and right-hand sides of \eqref{eq1a}. 
For $R$ we have: 
\begin{equation}
\label{eq1R}
R = \sum_{i=1}^k (k+j)[2(n+i) + k + j] = 
k(j + k)(j + 2k + 2n + 1).
\end{equation}
In particular, $R = k(k+1)$  if $j=0$. 

To simplify $L$, we use the well-known formula 
$1^2+ \dots +t^2 = (1/6)t(t+1)(2t+1)$  
for the sum of squares. 
We distinguish two cases: $m > n$ and $m \leq n$. 

\begin{equation}
\label{eq1L>} 
L_{m>n} = 
\frac{1}{6}\big[n(n+1)(2n+1) + (m-n)(m-n+1)(2(m-n)+1\big];  
\end{equation}

\begin{equation}
\label{eq1=<} 
L_{m \leq n} = 
\frac{1}{6}\big[(n(n+1)(2n+1) - (n-m-1)(n-m)(2(n-m)-1)\big]. 
\end{equation}

However, equation $L = R$ results in the same quadratic equation:

\begin{equation}
\label{eq2a} 
(m + 1)n^2 - [(2k(j+k) + m(m+1)]n - 
[k(j+k)(j+2k+1) - \frac{1}{6}m(2m^2+3m+1)] = 0,  
\end{equation}
whose discriminant $D = D(j,m,k)$ is given by formula 

\begin{equation}
\label{e-disc} 
D = 
4k\left(j(k + m + 1)(j + 2k + m + 1) + k(k + m + 1)^2 \right) - 
\frac{1}{3}m(m + 1)^2(m + 2). 
\end{equation}

\subsection{Case $m=0$} 
\subsubsection{Simplifications} 
In this case equations (\ref{eq2a}) and (\ref{e-disc}) 
are reduced to 
\begin{equation}
\label{eq2a-m=0} 
n^2 - 2k(j+k)n - k(j+k)(j+2k+1) = 0, \;\; 
D = D(j,0,k) = 4k(k+1)(k+j)(k+j+1). 
\end{equation}

The product $\Pi = \Pi(k,j) = k(k+1)(k+j)(k+j+1)$ 
for $k \geq 0$ and $j \geq 0$ is a quadratic number. 
By definition, it is the product of two triangular numbers and 4. 
Number $\Pi$ and the corresponding pair $(k,j)$ are called {\em square} if $\Pi$ is a perfect square. 

If $m=0$, equation (\ref{eq2a}) has two roots $n'$ and $n''$ given by formula   
\begin{equation}
\label{e-roots-m=0} 
n = n(k,j) = k(k+j) \pm \sqrt{\Pi} = 
\sqrt{k(k+j)}\big(\sqrt{k(k+j)} \pm \sqrt{(k+1)(k+j+1)}\big). 
\end{equation}

For any integer $k \geq 1$ and $j \geq 0$ we obtain the desired identity (\ref{eq1}) 
with $n = n(k,j)$ and $m=0$.
We are interested only in integer $n$
(which may be positive, negative, or 0).

\begin{lemma}
The following four statements are equivalent: 
\newline 
(i) one root of (\ref{eq2a-m=0})  is integer; \; 
(ii) both roots of (\ref{eq2a-m=0})  are integer;
\newline
(iii) discriminant $D(j,0,k) = 4\Pi(k,j)$ 
is a perfect square; \; 
\newline 
(iv) quadratic number $\Pi(k,j)$  
is a perfect square,  
that is, pair $(k,j)$  is square.
\end{lemma}

\proof It is straightforward. \qed 

\medskip 
 
If $j=0$, we simplify further as follows: 
$\Pi = k(k+1)(k+j)(k+j+1) = [k(k+1)]^2$,  
which is a perfect square for all $k$.   
Thus, we obtain Dostor's identity.  

Several other infinite families of identities
will be given in this section.

\subsubsection{Subcase ($m=0, \, k=1$)  and  
near-isosceles Pythagorean triples 
\newline 
$n^2 + (n+1)^2 = (n+j+2)^2$}
\label{near-iso}
The first such triples are 

\medskip 

(3, 4, 5), (20, 21, 29), (119, 120, 169), 
(696, 697, 985), (4,059, 4,060, 5,741), 

\smallskip 

(23,660, 23,661, 33,461), (137,903, 137,904, 195,025), 
(803,760, 803,761, 1,136,689),..  

\medskip 

In this case, (\ref{eq1}) is reduced to   
\begin{equation} 
\label{eq-m=0-k=1} 
n^2 = (j+1)(2n+j+3)\;\; 
\text{or, equivalently,} \;\;  n^2 = 2(j+1)n + (j+1)(j+3),
\end{equation}
which discriminant is  $D(j,0,1) = 8(j+1)(j+2)$ and solutions  
\begin{equation}
n=(2(j+1)\pm\sqrt{D})/2=(j+1)\pm\sqrt{2(j+1)(j+2)})
\end{equation}
define the near-isosceles Pythagorean triples $(n, n+1, n+j+2)$.

\medskip

Recall the``Fibonacci-like'' sequence 
$$(x_0,x_1, x_2, \dots) = 
(0, 1, 2, 5, 12, 29, 70, 169, 408, 985, \; 2,378, \; 5,741, \dots)$$   
initialized by 
$x_0 =  0, x_1 = 1$  and defined by recursion 
$x_{i+2} = 2 x_{i+1} + x_i$  for $i = 1, 2, \ldots$.
Lemma \ref{t-tr-sq} results in the following characterization. 

\begin{theorem}
\label{t-Pyth}
The $i$th near-isosceles Pythagorean triple 
\((n_i, n_i + 1, n_i + j_i + 2)\), for \(i = 1, 2, \dots\), 
is given by the following formulas:

\medskip

\[
n_i = x_{i+1}^2 - x_i^2 \quad \text{if } i \text{ is odd}, \qquad
n_i = 2x_i x_{i+1} \quad \text{if } i \text{ is even};
\]

\[
n_i + 1 = 2x_i x_{i+1} \quad \text{if } i \text{ is odd}, \qquad
n_i + 1 = x_{i+1}^2 - x_i^2 \quad \text{if } i \text{ is even};
\]

\[
j_i + 2 = 2x_i^2 \quad \text{if } i \text{ is odd}, \qquad
j_i + 2 = (x_{i+1} - x_i)^2 \quad \text{if } i \text{ is even}.
\]

\medskip

In both cases,
$\;\; n_i + j_i + 2 = x_i^2 + x_{i+1}^2, \;\;$
that is,
$\;\; j_i + 2 = x_i^2 + x_{i+1}^2 - n_i.$
\end{theorem}

The second solution $n'_i = (j_i+1) - \sqrt{2(j_i+1)(j_i+2)})$  
also provides the Pythagorean triple  
$(-n_{i-1}+1, -n_{i-1}, -(n_{i-1} + j_{i-1}+2))$ for $i>1$ 
and (-1, 0, 1) for $i=1$. 

\bigskip

The several first near-isosceles Pythagorean triples are shown in Table \ref{table-isosceles}.

\begin{table}[h] 
\small
\begin{center}
\begin{tabular}{ c c c c}
$j+1$ & $\sqrt{(j+1)(j+2)/2}$ & $n$ & $n^2+(n+1)^2=(n+j+2)^2$\\ 
\hline
$(0+1)^2 1^2 = 1$ & $1$ & $3$ & $3^2 + 4^2 = 5^2$\\
$2 \cdot 2^2 = 8$ & $6$ & $20$ & $20^2 + 21^2 = 29^2$\\
$(2+5)^2 = 7^2 =49$ & $35$ & $119$ & $119^2 + 120^2 = 169^2$\\
$2 \cdot 12^2 = 288$ & $204$ & $696$ & $696^2 + 697^2 = 985^2$\\
$(12+29)^2 = 41^2 = 1,681$ & $1,189$ & $4,059$ & $4,059^2 + 4,060^2 = 5,741^2$\\
$2 \cdot 70^2 = 9,800$ & $6,930$ & $23,660$ & $23,660^2 + 23,661^2 = 33,461^2$\\
$(70+169)^2 = 239^2 = 57,121$ & $40,391$ & $137,903$ & $137,903^2 + 137,904^2 = 195,025^2$\\
$2 \cdot 408^2 = 332,928$ & $235,416$ & $803,760$ & $803,760^2 + 803,761^2 = 1,136,689^2$\\
\hline
\end{tabular}
\end{center}
\caption{Several first near-isosceles Pythagorean triples} 
\label{table-isosceles}

\end{table}

\begin{remark}
In general the Pythagorean triplets 
$a^2 + b^2 = c^2$  are characterized by formulae  
$a^2 = y^2-z^2$, $b = 2yz$, $c = y^2+z^2$; 
see, for example, \cite{Hin83} for more details.
\end{remark}

\subsubsection{Bi-Pythagorean identities} 
In case  $m=0$, our main equation (\ref{eq1}) 
is reduced to \eqref{eq2a-m=0}. 
The latter has integer roots whenever 
its discriminant $D(j,0,k) = 4\Pi(k,j)$  is a perfect square. 

\begin{theorem}
\label{bi-Pyth}
Quadratic number  $\Pi(k,j) = k(k+1)(k+j)(k+j+1)$ 
is a perfect square whenever  
both triangular numbers 
$k(k+1)/2$ and $(k+j)(k+j+1)/2$ are perfect squares.   
\end{theorem} 

Such solutions of~\eqref{eq1} are called \emph{bi-Pythagorean}.
They are in a one-to-one correspondence with the square triangular numbers. 

\medskip

For \(k = 8\) and \(j = 41\), we have:
\[
\frac{k(k+1)}{2} = \frac{8 \cdot 9}{2} = 36 = 6^2 
\quad \text{and} \quad 
\frac{(k+j)(k+j+1)}{2} = \frac{49 \cdot 50}{2} = 1225 = 35^2.
\]

\medskip

Furthermore, 
\[
n = k(k + j) \pm \sqrt{\Pi} = 8 \cdot 49 \pm 2 \cdot 6 \cdot 35 = 392 \pm 420.
\]
So \(n' = 812\), \quad \(n'' = -28\), and equation~\eqref{eq1} yields the identities:
\[
812^2 + 813^2 + 814^2 + 815^2 + 816^2 + 817^2 + 818^2 + 819^2 + 820^2 = 
\]
\[
862^2 + 863^2 + 864^2 + 865^2 + 866^2 + 867^2 + 868^2 + 869^2 = 5,\!992,\!764.
\]

\[
(-28)^2 + (-27)^2 + (-26)^2 + (-25)^2 + (-24)^2 + (-23)^2 + (-22)^2 + (-21)^2 + (-20)^2 = 
\]
\[
22^2 + 23^2 + 24^2 + 25^2 + 26^2 + 27^2 + 28^2 + 29^2 = 5,\!244.
\]

The proof of the theorem is straightforward. 
Bi-Pythagorean identities are characterized 
by Theorem \ref{t-Pyth}.

\subsubsection{Main family of identities in case $m=0$} 
In this case we want to characterize all pairs $(k,j)$ 
such that the product   
$$\Pi(k,j) = k(k+1)(k+j)(k+j+1) = 
4 \frac{k(k+1)}{2} \frac{(k+j)(k+j+1)}{2},$$ 
that is, the product of two triangular numbers 
(times 4) is a perfect square. 
Obviously, if both numbers are square then the product also is. 
Yet, there are many other square quadratic numbers. 
For each of them we compute two roots   
$$n = k(k+j) \pm \sqrt{\Pi} = 
\sqrt{k(k+j)}\left[\sqrt{k(k+j)} \pm \sqrt{(k+1)(k+j+1)}\right]$$ 
and generate the corresponding quadratic identities. 
As we know, each of them provides two pairs of intervals  
of successive integers with equal sums of squares. 
Consider examples. 

\medskip 

For $t=1$ we have $j = j_1(k) = k(4k+3)$. 
Furthermore, $$k+j = 4k(k+1), \; k+j+1 = (2k+1)^2; \; 
\Pi_1(k) = k(k+1)(k+j)(k+j+1) = [2k(k+1)(2k+1)]^2.$$  
is a perfect square for all integer $k$ with  
$\sqrt{\Pi(k)} = 2k(k+1)(2k+1)$. 
Two roots of (\ref{eq2a-m=0}) are     
$$n = k(k+j) \pm \sqrt{\Pi} = 2k(k+1) [2k \pm (2k+1)], 
\;\; \text{that is}, \;\; 
n' = 2k(k+1)(4k+1), n'' = -2k(k+1).$$ 

The corresponding two identities are 
$$\sum_{t=0}^k [2k(k+1)(4k+1) + t]^2 = 
\sum_{t=1}^k [2k(k+1)(4k+3) + t]^2;$$
$$\sum_{t=0}^k [-2k(k+1) + t]^2 = 
\sum_{t=1}^k [2k(k+1) + t]^2.$$

For  $k = 1,2,3,4$ we obtain, respectively: 
$$20^2 + 21^2 = 29^2 = 841; \;\; 
(-4)^2 + (-3)^2 = 5^2 = 25;$$
$$108^2 + 109^2 + 110^2 = 133^2 + 134^2 = 35,645; \;\; 
(-12)^2 + (-11)^2 + (-10)^2 = 13^2 + 14^2 = 365;$$ 
$$312^2 + 313^2 + 314^2 + 315^2 = 361^2 + 362^2 + 363^2 = 393,134;$$
$$(-24)^2 + (-23)^2 + (-22)^2 + (-21)^2 = 25^2 + 26^2 + 27^2 = 2,030;$$ 
$$680^2 + 681^2 + 682^2 + 683^2 + 684^2 = 761^2 + 762^2 + 763^2 + 764^2 = 2,325,690;$$ 
$$(-40)^2 + (-39)^2 + (-38)^2 + (-37)^2 + (-36)^2 = 
41^2 + 42^2 + 43^2 + 44^2 = 7,230.$$

For positive integer $k$ and $i \leq 5$  we obtain similar identities:

\medskip 

For \( j = j_2(k) = k[(4k)^2 + 3(4k) + 8] = 8k(k+1)(2k+1) \) we have:

\smallskip 

\[
k + j = k(4k+3)^2 \quad \text{and} \quad k + j + 1 = (k+1)(4k+1)^2.
\]
Hence, the quadratic number \(\Pi_2(k)\) is a perfect square for all \(k\), with
\[
\sqrt{\Pi_2} = k(k+1)(4k+1)(4k+3).
\]
Two roots of \eqref{eq2a-m=0} are
\[
n = k(k+j) \pm \sqrt{\Pi} = k(4k+3) \bigl[k(4k+3) \pm (k+1)(4k+1)\bigr],
\]
that is,
\[
n' = k(4k+3)(8k^2 + 8k + 1) \quad \text{and} \quad n'' = -k(2k+1)(4k+3).
\]

\bigskip

For \( j = j_3(k) = k(4k+3)(16k^2 + 20k + 5) \) we have:

\[
k + j = 16k(k+1)(2k+1)^2 = 8 \sqrt{\Pi_1}, \quad k + j + 1 = (8k^2 + 8k + 1)^2.
\]

Hence,
\[
\sqrt{\Pi_3} = 4k(k+1)(2k+1)(8k^2 + 8k + 1),
\]
and two roots of \eqref{eq2a-m=0} are
\[
4k(k+1)(2k+1) \bigl[4k(2k+1) \pm (8k^2 + 8k + 1)\bigr],
\]
that is,
\[
n' = 4k(k+1)(2k+1)(16k^2 + 12k + 1), \quad n'' = -4k(k+1)(2k+1)(4k+1).
\]

\bigskip

For \( j = j_4(k) = 8k(k+1)(2k+1)(16k^2 + 16k + 3) \) we have:

\[
k + j = k \bigl[(4k)^2 + 5(4k) + 5\bigr]^2, \quad
k + j + 1 = (k+1) \bigl[(4k)^2 + 3(4k) + 1\bigr]^2.
\]

Hence, \(\Pi = \Pi(k, j_4(k))\) is a perfect square for all \(k\), with
\[
\sqrt{\Pi_4} = k(k+1)(16k^2 + 12k + 1)(16k^2 + 20k + 5).
\]
Two roots of \eqref{eq2a-m=0} are
\[
n = k(16k^2 + 20k + 5) \bigl[k(16k^2 + 20k + 5) \pm (k+1)(16k^2 + 12k + 1)\bigr],
\]
that is,
\[
n' = k(2k+1)(16k^2 + 16k + 1)(16k^2 + 20k + 5),
\quad
n'' = -k(8k^2 + 8k + 1)(16k^2 + 20k + 5).
\]

\medskip

For \( j = j_5(k) = k(16k^2 + 20k + 5)(64k^3 + 112k^2 + 56k + 7) \) we have

\[
k + j = 4k(k+1)(4k+1)^2 (4k+3)^2, \quad
k + j + 1 = (2k+1)^2 (16k^2 + 16k + 1)^2.
\]

Hence, \(\Pi_5(k)\) is a perfect square for all integer \(k\), with
\[
\sqrt{\Pi_5} = 2k(k+1)(2k+1)(4k+1)(4k+3)(16k^2 + 16k + 1).
\]
Two roots of \eqref{eq2a-m=0} are
\[
n' = 2k(k+1)(4k+1)(4k+3)(64k^3 + 80k^2 + 24k + 1),
\quad
n'' = -2k(k+1)(4k+1)(4k+3)(16k^2 + 12k + 1).
\]

\subsection{Case $k=1$}
In general,  without the assumption that $m=0$, this case is reduced to the equation 

\smallskip 
\begin{equation}
(n-m)^2 + (n-m+1)^2 + \dots + (n+1)^2 = (n+j+2)^2,
\label{Lucas}
\end{equation}

\smallskip 
\noindent 
whose right-hand side contains only one square. 
Solving this equation is 
``intimately related to finding 
all integral points on elliptic curves belonging to a certain family which can be represented by a
Weierstrass equation with parameter $n-m$'' \cite{BST97}.

The complete solution of \eqref{Lucas} is presented in~\cite{BST97}.  
The special case $m = n$ is known as the \emph{Lucas' Square Pyramid Problem}.  
For example, the following identity holds:
\[
1^2 + 2^2 + \dots + 24^2 = 70^2.
\]
This subcase has attracted considerable attention over the years; 
the reader can find  a comprehensive list of historical and contemporary references, 
in \cite[ch. 6, Problem D3]{Guy94}. 

\subsection{Case $j=0$ and $m>0$} 
\label{s-m>0}
If $j=0$, the discriminant of (\ref{eq1}) is reduced to   
$$D(0,m,k) = [2k(k + m + 1)]^2 - \frac{1}{3}m(m + 1)^2(m + 2),$$  
that is,  $D$ is a perfect square minus a function of $m$. 
If also $m=0$ then  $D(0,0,k) = [2k(k+1)]^2$  is a perfect square 
and we obtain Theorem \ref{t-Dos}. 
Partial results for $m>0$ are given below, in Table \ref{t4}.
The following table represents all rational solutions
for $1 \leq m \leq 100, \; 1 \leq k \leq 100$  and $j=0$. 
As we know, the denominators are divisors of $2(m+1)$. 
We are interested only in integer solutions.

\begin{table}[h]
\small
\begin{center}
\begin{tabular}{ c c c c}
$m$ & $k$ & $n$ & $\sqrt{D}$\\ 
\hline
7 & 2 & 3 & 16\\
7 & 2 & 5 & 16\\
8 & 8 & -34/9  &268\\
8 & 8 & 26 & 268\\
17 & 7 & 3 & 296\\
17 & 7 & 175/9 & 296\\
24 & 10 & 4 & 600\\
24 & 10 & 28 & 600\\
25 &13 & 1 &936\\
25 & 13 & 37 &936\\
25 & 25 & -155/13 &2,520\\
25 & 25 & 85 &2,520\\
53&  39  &-32/3 &7,056\\
53&  39  & 120&7,056\\
55&  17&  35/2 &1,698\\
55&  17&  1339/28 &1,698\\
62&  22&  137/9 &2,956\\
62&  22&  435/7 &2,956\\
89&  37&  649/45 &8,152\\
89&  37&  105 &8,152\\

\hline
\end{tabular}
\end{center}
\caption{All rational solutions of \eqref{eq1} for $j=0$, $1 \leq m \leq 100, \; 1 \leq k \leq 100$. }  \label{t4}
\end{table}

It is not difficult to verify that 
inequalities  $0 < m < 7$  cannot hold for a solution.
Indeed, the determinant $D(0,m,k)$ given by (\ref{e-disc}) 
is a complete square minus $\frac{1}{3}m(m+1)^2(m+2)$. 
But $D(0,m,k)$ itself must be a complete square. 
This implies an upper bound for $k$.
In particular, for  $m \leq 5$  we obtain: 

\smallskip 

$D(0,0,k) = (2k(k+1))^2$,  

\smallskip 

$D(1,0,k) = (2k(k+2))^2 - 4$, hence, $k(k+2) <  1 \frac{1}{4}$,  no such $k$; 

\smallskip 

$D(2,0,k) = (2k(k+3))^2 - 24, \; k(k+3) < 6 \frac{1}{4}$, hence,  $k = 1$; 

\smallskip 

$D(3,0,k) = (2k(k+4))^2 - 80,  \; k(k+4) < 20 \frac{1}{4}$, hence, $k \leq 2$;    
\smallskip 

$D(4,0,k) = (2k(k+4))^2 - 200, \; k(k+5) < 50 \frac{1}{4}$, hence, $k \leq 5$;  

\smallskip 

$D(5,0,k) = (2k(k+4))^2 - 420,  \;  k(k+6) < 105 \frac{1}{4}$, hence, $k \leq 7$;

\smallskip 

$D(6,0,k) = (2k(k+4))^2 - 784, \; k(k+7) < 196 \frac{1}{4}$, hence, $k \leq 10$. 

\medskip 

On the other hand, $m$  cannot be much larger than  $k$. 
Indeed, the discriminant must be nonnegative,  
$D(0,m,k) \geq 0$, and hence, inequality
$k(k+m+1) \geq \frac{m+1}{2} \sqrt{\frac{m(m+2)}{3}}$  holds.
For example, for $m=7$ and $k=2$  we have, $20 > 4 \sqrt{21}$.

\subsection{Data for $m \leq 100, k \leq 100$ and $j \leq 18$} 
We collect all cases for $1 \leq m \leq 100, 1 \leq k \leq 100$ and $1 \leq j \leq 18$ when discriminant $D(j,m,k)$ is a perfect square (see \ref{t5}). 
\vspace{-0.8em}
\begin{table}[h]
\footnotesize
\centering
\begin{minipage}[t]{0.32\textwidth}
\centering
\caption*{\footnotesize $j = 1$}
\vspace{-0.8em}
\begin{tabular}{|c|c|c|c|}
\hline
$m$ & $k$ & $n$ & $\sqrt{D}$ \\
\hline
7 & 2 & $2$ & $36$ \\
7 & 2 & $13/2$ & $36$ \\
16 & 4 & $108/17$ & $96$ \\
16 & 4 & $12$ & $96$ \\
22 & 22 & $-11$ & $2,024$ \\
22 & 22 & $77$ & $2,024$ \\
31 & 7 & $17$ & $16$ \\
31 & 7 & $35/2$ & $16$ \\
31 & 13 & $35/8$ & $1,076$ \\
31 & 13 & $38$ & $1,076$ \\
31 & 24 & $-8$ & $2,704$ \\
31 & 24 & $153/2$ & $2,704$ \\
54 & 32 & $-3$ & $5,412$ \\
54 & 32 & $477/5$ & $5,412$ \\
62 & 15 & $85/3$ & $816$ \\
62 & 15 & $289/7$ & $816$ \\
79 & 35 & $19/2$ & $7,320$ \\
79 & 35 & $101$ & $7,320$ \\
\hline
\end{tabular}
\end{minipage}
\hfill
\begin{minipage}[t]{0.32\textwidth}
\centering
\caption*{\footnotesize $j = 2$}
\vspace{-0.8em}
\begin{tabular}{|c|c|c|c|}
\hline
$m$ & $k$ & $n$ & $\sqrt{D}$ \\
\hline
5 & 5 & $-10/3$ & $140$ \\
5 & 5 & $20$ & $140$ \\
9 & 5 & $-1$ & $180$ \\
9 & 5 & $17$ & $180$ \\
21 & 9 & $2$ & $572$ \\
21 & 9 & $28$ & $572$ \\
22 & 8 & $91/23$ & $484$ \\
22 & 8 & $25$ & $484$ \\
31 & 38 & $-23$ & $5,504$ \\
31 & 38 & $149$ & $5,504$ \\
40 & 15 & $305/41$ & $1,540$ \\
40 & 15 & $45$ & $1,540$ \\
48 & 49 & $-25$ & $9,800$ \\
48 & 49 & $175$ & $9,800$ \\
62 & 14 & $31$ & $448$ \\
62 & 14 & $343/9$ & $448$ \\
62 & 15 & $80/3$ & $1,056$ \\
62 & 15 & $304/7$ & $1,056$ \\
63 & 14 & $35$ & $0$ \\
\hline
\end{tabular}
\end{minipage}
\hfill
\begin{minipage}[t]{0.32\textwidth}
\centering
\caption*{\footnotesize $j = 3$}
\vspace{-0.8em}
\begin{tabular}{|c|c|c|c|}
\hline
$m$ & $k$ & $n$ & $\sqrt{D}$ \\
\hline
16 & 11 & $-66/17$ & $712$ \\
16 & 11 & $38$ & $712$ \\
24 & 25 & $-14$ & $2,700$ \\
24 & 25 & $94$ & $2,700$ \\
30 & 27 & $-395/31$ & $3,340$ \\
30 & 27 & $95$ & $3,340$ \\
31 & 10 & $57/8$ & $796$ \\
31 & 10 & $32$ & $796$ \\
32 & 10 & $260/33$ & $796$ \\
32 & 10 & $32$ & $796$ \\
47 & 83 & $-727/12$ & $22,348$ \\
47 & 83 & $405$ & $22,348$ \\
55 & 14 & $41/2$ & $1,260$ \\
55 & 14 & $43$ & $1,260$ \\
62 & 14 & $257/9$ & $784$ \\
62 & 14 & $41$ & $784$ \\
62 & 25 & $83/9$ & $4,144$ \\
62 & 25 & $75$ & $4,144$ \\
62 & 71 & $-283/7$ & $19,508$ \\
62 & 71 & $2423/9$ & $19,508$ \\
70 & 35 & $175/71$ & $7,280$ \\
70 & 35 & $105$ & $7,280$ \\
71 & 71 & $-71/2$ & $20,732$ \\
71 & 71 & $2272/9$ & $20,732$ \\
79 & 36 & $71/10$ & $7,992$ \\
79 & 36 & $107$ & $7,992$ \\
\hline
\end{tabular}
\end{minipage}
\end{table}
\vspace{-0.9em}
\begin{table}[!htbp]
\footnotesize
\centering
\begin{minipage}[t]{0.32\textwidth}
\centering
\caption*{\footnotesize $j = 4$}
\vspace{-0.8em}
\begin{tabular}{|c|c|c|c|}
\hline
$m$ & $k$ & $n$ & $\sqrt{D}$ \\
\hline
7 & 7 & $-21/4$ & $294$ \\
7 & 7 & $63/2$ & $294$ \\
9 & 1 & $5$ & $0$ \\
17 & 3 & $7$ & $96$ \\
17 & 3 & $37/3$ & $96$ \\
31 & 6 & $29/2$ & $184$ \\
31 & 6 & $81/4$ & $184$ \\
31 & 10 & $13/2$ & $856$ \\
31 & 10 & $133/4$ & $856$ \\
31 & 19 & $-67/16$ & $2,134$ \\
31 & 19 & $125/2$ & $2,134$ \\
45 & 69 & $-48$ & $16,560$ \\
45 & 69 & $312$ & $16,560$ \\
54 & 11 & $30$ & $0$ \\
\hline
\end{tabular}
\end{minipage}
\hfill
\begin{minipage}[t]{0.32\textwidth}
\centering
\caption*{\footnotesize $j = 5$}
\vspace{-0.8em}
\begin{tabular}{|c|c|c|c|}
\hline
$m$ & $k$ & $n$ & $\sqrt{D}$ \\
\hline
7 & 2 & $0$ & $84$ \\
7 & 2 & $21/2$ & $84$ \\
8 & 3 & $-2/3$ & $132$ \\
8 & 3 & $14$ & $132$ \\
23 & 14 & $-4$ & $1,276$ \\
23 & 14 & $295/6$ & $1,276$ \\
30 & 58 & $-1403/31$ & $11,044$ \\
30 & 58 & $311$ & $11,044$ \\
40 & 42 & $-972/41$ & $7,532$ \\
40 & 42 & $160$ & $7,532$ \\
48 & 14 & $90/7$ & $1,624$ \\
48 & 14 & $46$ & $1,624$ \\
79 & 41 & $3/20$ & $10,068$ \\
79 & 41 & $126$ & $10,068$ \\
\hline
\end{tabular}
\end{minipage}
\hfill
\begin{minipage}[t]{0.32\textwidth}
\centering
\caption*{\footnotesize $j = 6$}
\vspace{-0.8em}
\begin{tabular}{|c|c|c|c|}
\hline
$m$ & $k$ & $n$ & $\sqrt{D}$ \\
\hline
9 & 1 & $3$ & $44$ \\
9 & 1 & $37/5$ & $44$ \\
21 & 3 & $10$ & $76$ \\
21 & 3 & $148/11$ & $76$ \\
24 & 4 & $10$ & $180$ \\
24 & 4 & $86/5$ & $180$ \\
32 & 32 & $-604/33$ & $4,696$ \\
32 & 32 & $124$ & $4,696$ \\
48 & 35 & $-87/7$ & $6,440$ \\
48 & 35 & $119$ & $6,440$ \\
54 & 16 & $69/5$ & $2,156$ \\
54 & 16 & $53$ & $2,156$ \\
57 & 11 & $29$ & $316$ \\
57 & 11 & $999/29$ & $316$ \\
79 & 84 & $-46$ & $28,800$ \\
79 & 84 & $314$ & $28,800$ \\
\hline
\end{tabular}
\end{minipage}
\end{table}

\vspace{-0.8em}
\begin{table}[h]
\footnotesize
\centering
\begin{minipage}[t]{0.32\textwidth}
\centering
\caption*{\footnotesize $j = 7$}
\vspace{-0.8em}
\begin{tabular}{|c|c|c|c|}
\hline
$m$ & $k$ & $n$ & $\sqrt{D}$ \\
\hline
22 & 23 & $-15$ & $2,576$ \\
22 & 23 & $97$ & $2,576$ \\
31 & 17 & $-7/2$ & $2,032$ \\
31 & 17 & $60$ & $2,032$ \\
32 & 25 & $-380/33$ & $3,416$ \\
32 & 25 & $92$ & $3,416$ \\
48 & 77 & $-56$ & $20,776$ \\
48 & 77 & $368$ & $20,776$ \\
70 & 94 & $-4371/71$ & $32,700$ \\
70 & 94 & $399$ & $32,700$ \\
71 & 15 & $92/3$ & $1,356$ \\
71 & 15 & $99/2$ & $1,356$ \\
71 & 22 & $155/9$ & $3,908$ \\
71 & 22 & $143/2$ & $3,908$ \\
\hline
\end{tabular}
\end{minipage}
\hfill
\begin{minipage}[t]{0.32\textwidth}
\centering
\caption*{\footnotesize $j = 8$}
\vspace{-0.8em}
\begin{tabular}{|c|c|c|c|}
\hline
$m$ & $k$ & $n$ & $\sqrt{D}$ \\
\hline
7 & 3 & $-5/2$ & $162$ \\
7 & 3 & $71/4$ & $162$ \\
21 & 3 & $8$ & $176$ \\
21 & 3 & $16$ & $176$ \\
31 & 6 & $21/2$ & $488$ \\
31 & 6 & $103/4$ & $488$ \\
32 & 5 & $157/11$ & $244$ \\
32 & 5 & $65/3$ & $244$ \\
32 & 11 & $11/3$ & $1,232$ \\
32 & 11 & $41$ & $1,232$ \\
41 & 17 & $47/21$ & $2,384$ \\
41 & 17 & $59$ & $2,384$ \\
49 & 27 & $-21/5$ & $4,760$ \\
49 & 27 & $91$ & $4,760$ \\
56 & 38 & $-12$ & $8,056$ \\
56 & 38 & $388/3$ & $8,056$ \\
62 & 20 & $115/9$ & $3,416$ \\
62 & 20 & $67$ & $3,416$ \\
85 & 43 & $0$ & $11,696$ \\
85 & 43 & $136$ & $11,696$ \\
94 & 49 & $-6/5$ & $14,744$ \\
94 & 49 & $154$ & $14,744$ \\
\hline
\end{tabular}
\end{minipage}
\hfill
\begin{minipage}[t]{0.32\textwidth}
\centering
\caption*{\footnotesize $j = 9$}
\vspace{-0.8em}
\begin{tabular}{|c|c|c|c|}
\hline
$m$ & $k$ & $n$ & $\sqrt{D}$ \\
\hline
7 & 2 & $-3/2$ & $124$ \\
7 & 2 & $14$ & $124$ \\
16 & 5 & $4/17$ & $404$ \\
16 & 5 & $24$ & $404$ \\
16 & 53 & $-840/17$ & $8,524$ \\
16 & 53 & $452$ & $8,524$ \\
23 & 26 & $-37/2$ & $3,260$ \\
23 & 26 & $352/3$ & $3,260$ \\
24 & 10 & $-4/5$ & $1,020$ \\
24 & 10 & $40$ & $1,020$ \\
31 & 24 & $-12$ & $3,344$ \\
31 & 24 & $185/2$ & $3,344$ \\
31 & 40 & $-57/2$ & $6,736$ \\
31 & 40 & $182$ & $6,736$ \\
54 & 24 & $9/5$ & $4,356$ \\
54 & 24 & $81$ & $4,356$ \\
79 & 30 & $11$ & $6,900$ \\
79 & 30 & $389/4$ & $6,900$ \\
79 & 51 & $-13$ & $14,520$ \\
79 & 51 & $337/2$ & $14,520$ \\
\hline
\end{tabular}
\end{minipage}
\end{table}
\vspace{-0.9em}
\begin{table}[!htbp]
\footnotesize
\centering
\begin{minipage}[t]{0.32\textwidth}
\centering
\caption*{\footnotesize $j = 10$}
\vspace{-0.8em}
\begin{tabular}{|c|c|c|c|}
\hline
$m$ & $k$ & $n$ & $\sqrt{D}$ \\
\hline
25 & 3 & $11$ & $156$ \\
25 & 3 & $17$ & $156$ \\
31 & 8 & $6$ & $896$ \\
31 & 8 & $34$ & $896$ \\
47 & 76 & $-57$ & $20,800$ \\
47 & 76 & $1129/3$ & $20,800$ \\
48 & 11 & $101/7$ & $1,400$ \\
48 & 11 & $43$ & $1,400$ \\
48 & 18 & $32/7$ & $2,912$ \\
48 & 18 & $64$ & $2,912$ \\
48 & 49 & $-29$ & $10,976$ \\
48 & 49 & $195$ & $10,976$ \\
80 & 27 & $15$ & $6,048$ \\
80 & 27 & $269/3$ & $6,048$ \\
94 & 47 & $0$ & $14,288$ \\
94 & 47 & $752/5$ & $14,288$ \\
\hline
\end{tabular}
\end{minipage}
\hfill
\begin{minipage}[t]{0.32\textwidth}
\centering
\caption*{\footnotesize $j = 11$}
\vspace{-0.8em}
\begin{tabular}{|c|c|c|c|}
\hline
$m$ & $k$ & $n$ & $\sqrt{D}$ \\
\hline
23 & 14 & $-41/6$ & $1,580$ \\
23 & 14 & $59$ & $1,580$ \\
31 & 7 & $7$ & $796$ \\
31 & 7 & $255/8$ & $796$ \\
31 & 8 & $11/2$ & $944$ \\
31 & 8 & $35$ & $944$ \\
62 & 11 & $193/7$ & $916$ \\
62 & 11 & $379/9$ & $916$ \\
62 & 56 & $-29$ & $15,064$ \\
62 & 56 & $1891/9$ & $15,064$ \\
79 & 22 & $83/4$ & $4,452$ \\
79 & 22 & $382/5$ & $4,452$ \\
80 & 43 & $-14/3$ & $11,880$ \\
80 & 43 & $142$ & $11,880$ \\
\hline
\end{tabular}
\end{minipage}
\hfill
\begin{minipage}[t]{0.32\textwidth}
\centering
\caption*{\footnotesize $j = 12$}
\vspace{-0.8em}
\begin{tabular}{|c|c|c|c|}
\hline
$m$ & $k$ & $n$ & $\sqrt{D}$ \\
\hline
14 & 2 & $7/3$ & $196$ \\
14 & 2 & $77/5$ & $196$ \\
16 & 5 & $-1$ & $476$ \\
16 & 5 & $27$ & $476$ \\
23 & 2 & $37/3$ & $16$ \\
23 & 2 & $13$ & $16$ \\
25 & 13 & $-5$ & $1,560$ \\
25 & 13 & $55$ & $1,560$ \\
31 & 7 & $13/2$ & $842$ \\
31 & 7 & $525/16$ & $842$ \\
55 & 41 & $-71/4$ & $9,414$ \\
55 & 41 & $2105/14$ & $9,414$ \\
57 & 17 & $11$ & $3,016$ \\
57 & 17 & $63$ & $3,016$ \\
62 & 13 & $150/7$ & $1,856$ \\
62 & 13 & $458/9$ & $1,856$ \\
78 & 67 & $-32$ & $21,804$ \\
78 & 67 & $244$ & $21,804$ \\
80 & 30 & $10$ & $7,380$ \\
80 & 30 & $910/9$ & $7,380$ \\
95 & 19 & $1957/48$ & $2,470$ \\
95 & 19 & $133/2$ & $2,470$ \\
\hline
\end{tabular}
\end{minipage}
\end{table}

\vspace{-0.8em}
\begin{table}[h]
\footnotesize
\centering
\begin{minipage}[t]{0.32\textwidth}
\centering
\caption*{\footnotesize $j = 13$}
\vspace{-0.8em}
\begin{tabular}{|c|c|c|c|}
\hline
$m$ & $k$ & $n$ & $\sqrt{D}$ \\
\hline
16 & 4 & $0$ & $408$ \\
16 & 4 & $24$ & $408$ \\
54 & 9 & $111/5$ & $924$ \\
54 & 9 & $39$ & $924$ \\
70 & 13 & $2025/71$ & $1,596$ \\
70 & 13 & $51$ & $1,596$ \\
80 & 14 & $38$ & $1,080$ \\
80 & 14 & $154/3$ & $1,080$ \\
\hline
\end{tabular}
\end{minipage}
\hfill
\begin{minipage}[t]{0.32\textwidth}
\centering
\caption*{\footnotesize $j = 14$}
\vspace{-0.8em}
\begin{tabular}{|c|c|c|c|}
\hline
$m$ & $k$ & $n$ & $\sqrt{D}$ \\
\hline
7 & 3 & $-19/4$ & $234$ \\
7 & 3 & $49/2$ & $234$ \\
21 & 11 & $-6$ & $1,276$ \\
21 & 11 & $52$ & $1,276$ \\
31 & 3 & $267/16$ & $26$ \\
31 & 3 & $35/2$ & $26$ \\
38 & 77 & $-194/3$ & $20,540$ \\
38 & 77 & $462$ & $20,540$ \\
48 & 21 & $-1$ & $3,920$ \\
48 & 21 & $79$ & $3,920$ \\
62 & 16 & $15$ & $2,976$ \\
62 & 16 & $1307/21$ & $2,976$ \\
89 & 37 & $19/3$ & $10,644$ \\
89 & 37 & $623/5$ & $10,644$ \\
\hline
\end{tabular}
\end{minipage}
\hfill
\begin{minipage}[t]{0.32\textwidth}
\centering
\caption*{\footnotesize $j = 15$}
\vspace{-0.8em}
\begin{tabular}{|c|c|c|c|}
\hline
$m$ & $k$ & $n$ & $\sqrt{D}$ \\
\hline
22 & 31 & $-27$ & $4,600$ \\
22 & 31 & $173$ & $4,600$ \\
32 & 17 & $-20/3$ & $2,584$ \\
32 & 17 & $788/11$ & $2,584$ \\
48 & 6 & $176/7$ & $140$ \\
48 & 6 & $28$ & $140$ \\
48 & 83 & $-66$ & $25,088$ \\
48 & 83 & $446$ & $25,088$ \\
95 & 61 & $-209/12$ & $21,736$ \\
95 & 61 & $209$ & $21,736$ \\
\hline
\end{tabular}
\end{minipage}
\end{table}
\begin{table}[!htbp]
\footnotesize
\centering
\begin{minipage}[t]{0.32\textwidth}
\centering
\caption*{\footnotesize $j = 16$}
\vspace{-0.8em}
\begin{tabular}{|c|c|c|c|}
\hline
$m$ & $k$ & $n$ & $\sqrt{D}$ \\
\hline
7 & 7 & $-21/2$ & $546$ \\
7 & 7 & $231/4$ & $546$ \\
31 & 20 & $-11$ & $3,136$ \\
31 & 20 & $87$ & $3,136$ \\
62 & 14 & $17$ & $2,604$ \\
62 & 14 & $175/3$ & $2,604$ \\
62 & 26 & $1$ & $5,964$ \\
62 & 26 & $287/3$ & $5,964$ \\
79 & 69 & $-71/2$ & $23,730$ \\
79 & 69 & $2089/8$ & $23,730$ \\
\hline
\end{tabular}
\end{minipage}
\hfill
\begin{minipage}[t]{0.32\textwidth}
\centering
\caption*{\footnotesize $j = 17$}
\vspace{-0.8em}
\begin{tabular}{|c|c|c|c|}
\hline
$m$ & $k$ & $n$ & $\sqrt{D}$ \\
\hline
23 & 2 & $8$ & $244$ \\
23 & 2 & $109/6$ & $244$ \\
31 & 3 & $51/4$ & $296$ \\
31 & 3 & $22$ & $296$ \\
47 & 83 & $-135/2$ & $25,336$ \\
47 & 83 & $1381/3$ & $25,336$ \\
54 & 15 & $9$ & $2,940$ \\
54 & 15 & $687/11$ & $2,940$ \\
\hline
\end{tabular}
\end{minipage}
\hfill
\begin{minipage}[t]{0.32\textwidth}
\centering
\caption*{\footnotesize $j = 18$}
\vspace{-0.8em}
\begin{tabular}{|c|c|c|c|}
\hline
$m$ & $k$ & $n$ & $\sqrt{D}$ \\
\hline
14 & 2 & $1/3$ & $280$ \\
14 & 2 & $19$ & $280$ \\
57 & 11 & $17$ & $1,972$ \\
57 & 11 & $51$ & $1,972$ \\
62 & 8 & $2117/63$ & $88$ \\
62 & 8 & $35$ & $88$ \\
79 & 12 & $38$ & $960$ \\
79 & 12 & $50$ & $960$ \\
80 & 12 & $40$ & $720$ \\
80 & 12 & $440/9$ & $720$ \\
\hline
\end{tabular}
\end{minipage}
\caption{All cases when $D(j,m,k)$ is a perfect square for\\ $m = 1, \ldots, 100$, $k = 1, \ldots, 100$, $j = 1, \ldots, 18$.}
\label{t5}
\end{table}

\clearpage
\section*{Acknowledgements} 
The authors thank Rashit Ziatdinov 
(a chess grand-master from Tashkent, Uzbekistan)   
who brought our attention to the identity considered 
in case $m=j=0, k=6$ and to Besfort Shala for reference \cite{IMC25}. 
This paper is an output of the research project (HSE-BR-2025-024) 
implemented as part of the Basic Research Program at HSE University.

\end{document}